\documentstyle[11pt]{article}

\def\qed{\hfill $\Box$}
\def\rz{ I \!\! R}
\def\cz{ I \!\!\!\! C}

\def\gz{ Z \!\!\! Z}
\def\pr{ I \!\! P}

\def\p{\partial}
\def\s{\subset}
\def\se{\setminus}

\def\o{\overline}

\def\be{\beta}
\def\g{{\gamma}}
\def\d{{\delta}}
\def\De{{\Delta}}

\def\ph{\varphi}
\def\la{\lambda}

\def\la{\lambda}
\def\o{{\omega}}
\def\si{\sigma}
\def\Ga{\Gamma}

\def\M{{\cal M}}
\def\Z{{\cal Z}}
\def\D{{\cal D}}
\def\J{{\cal J}}
\def\S{{\cal S}}

\def\Si{{\Sigma}}
\def\noin{\noindent}
\def\be{\begin{equation}}
\def\ee{\end{equation}}
\def\vs{\bigskip}
\def\ms{\medskip}
\def\ss{\smallskip}

\begin{document}

\title {Plane curves with a big fundamental group of the  complement}

\vspace{.5cm}

\author{ G. Dethloff, S. Orevkov \footnote{Partually supported by Russian
Foundation for Fundamental
Research, N96-01-01218}, M. Zaidenberg}

\date{}
\maketitle

\begin{abstract}
\noin Let $C \s \pr^2$ be an irreducible plane curve
whose dual $C^* \s \pr^{2*}$ is an immersed curve
which is neither a conic nor a nodal cubic.
The main result states that the Poincar\'e group $G=\pi_1(\pr^2 \se C)$
contains a free group with two generators.
If the geometric genus $g$ of $C$ is at least $2$, then 
a subgroup of $G$ can be mapped epimorphically onto
the fundamental group of the normalization of $C$,
and the result follows.
To handle the cases $g=0,1$,
we construct universal families of
immersed plane curves and their Picard bundles. This allows us to reduce the
consideration to the case of Pl\"ucker curves. Such a curve $C$ can be regarded
as a plane section of the corresponding discriminant
hypersurface (cf. [Zar, DoLib]). Applying Zariski--Lefschetz type arguments
we deduce the result from `the bigness' of the braid group $B_{d,\,g}$, that
is, of the group of $d$--string braids of a compact genus $g$ Riemann surface.

\end{abstract}

\section*{Introduction}

\vs

\noin The fundamental groups of the plane curve complements are of permanent
interest (see e.g. [Di, DoLib,  Lib,  MoTe,  No, O,  Zar] and the literature
therein). Here we look for the most coarse properties of these groups (cf. e.g.
[MoTe]). Namely, we distinguish between {\it big} and {\it small} groups.

\vs

\noin {\bf 0.1. Definition.} We say that a group $G$ is {\it big} if it
contains a non--abelian free subgroup. We call $G$ {\it small} if it
is {\it almost solvable}, i.e. it has a solvable subgroup of finite index.

\vs

Recall the Tits alternative [Ti]: {\it any subgroup $G$ of a general linear
group $GL(n,\,k)$ over a field $k$ of characteristic zero is either
big or small.} This alternative holds true, even in a stronger form, for  some
classes of discrete groups, such as hyperbolic groups in
sense of Gromov and the mapping class groups (see sect.1 below for references).

\ss

In [MoTe] classes of plane Pl\"ucker curves were indicated with
infinite almost solvable (i.e. small) non--abelian fundamental groups
of the complement. An example
is the branching divisor of a generic projection of the Veronese surface
$V_3$ of order $3$ onto $\pr^2$ [MoTe].

The well known Deligne--Fulton Theorem asserts that the complement of a nodal
plane curve has abelian fundamental group. Here we show (and this is the main
purpose of the paper) that
the fundamental group of
the complement of the dual of a nodal plane curve is big
(with two evident exceptions). More precisely, we have

\vs

\noin {\bf 0.2. Theorem.} {\it  Let $C \s \pr^2$ be an irreducible immersed
curve\footnote{i.e. all the analytic branches at the singular points of $C$
are smooth.}
which is neither a line nor a conic nor a nodal cubic.
Let $C^* \s \pr^{2*}$ be the dual curve.
Then the group $\pi_1( \pr^{2*} \se C^*)$ is big.}

\vs

\noin Obviously, the statement does not hold for
a line, nor for a conic.
If $C$ is a nodal cubic then
$C^*$ is a three--cuspidal quartic and $\pi_1 (\pr^{2*} \se C^*)$ is the
metacyclic group of order $12$ [Zar, p.143--145].

Let $d$ be the degree and $g$ be the geometric genus of $C$.
For $g\ge2$ the proof of Theorem 0.2 is easy.
Indeed,
denote by $\nu: C_{\rm norm} \to C$ a normalization of $C$.
Set
$R = \{(p,l)\in C_{\rm norm} \times \pr^{2*}\,|\,\nu(p)\in l\}$.
It is a smooth surface. Let $\mu_1 : R \to C_{\rm norm}$
and $\mu_2 : R \to \pr^{2*}$ be the canonical
projections. Since $C$ is immersed, it is easily seen that
$\mu_2$ is ramified exactly over $C^*$. Denote by
$R_0$ the part of $R$ over $\pr^{2*} \setminus C^*$.
Then $\mu_1 : R_0 \to C_{\rm norm}$ is a holomorphic surjection
with connected fibres. It follows that
$(\mu_1)_* : \pi_1(R_0) \to \pi_1(C_{\rm norm})$ is an epimorphism,
and hence the group $\pi_1(R_0)$ is big as soon as
$g(C_{\rm norm}) \ge 2$ (see e.g. 1.2(a) below).
Since $\mu_2 : R_0 \to \pr^{2*} \setminus C^*$
is a finite unramified covering, we have that
$(\mu_2)_* (\pi_1(R_0))$ is a finite index subgroup of
$\pi_1(\pr^{2*} \setminus C^*)$.
Therefore, the group $\pi_1(\pr^{2*} \setminus C^*)$ is also big.

Thus, the only non-trivial cases are $g=0$ and $g=1$.
However, the proofs of most of the intermediate results
needed for these two cases
are valid for any $g$, some of them under the additional assumption that
$d\ge 2g-1$ (which is automatically fulfilled for $g = 0,\,1$).
Therefore, we formulate everything for an arbitrary genus.
This provides another proof of Theorem 0.2
for the case $g\ge2$, $d\ge 2g-1$.

\ss

The paper is organized as follows.
In Section 1 we provide some (mostly well known)
examples of big groups.
Besides, by several examples we illustrate a conjectural
relation between bigness of the fundamental group and C-hyperbolicity.
These include, in particular,
the quasi--projective
quotients of bounded symmetric domains and the complements
of certain reducible plane curves.

The proof of Theorem 0.2 is done in Sect. 4. The
results in Sect. 2 and 3 (which we believe to be of some independent interest)
reduce the proof to the
case of a nodal Pl\"ucker curve. In Theorem 2.1 we show that the part
$Imm_{d,\,g}$ of
the Hilbert scheme of degree $d$ genus $g$ plane curves, which corresponds to
the immersed curves, is smooth,
and the universal family of curves admits a simultaneous normalization
over $Imm_{d,\,g}$ (see [AC, Ha] for related results, especially concerning
($a$) and
($b$) of Theorem 2.1). We show also that the nodal and
(for $d \ge 2g-1$) the Pl\"ucker curves form Zariski open subsets of
$Imm_{d,\,g}$.

\ms

A preliminary version of Theorem 0.2
(under the additional restriction $d\ge 2g-1$)
was announced in [DeZa1, Sect. 7] (see also [DeZa2]).
%Note that a presentation of the group $\pi_1(\pr^{2*} \se C^*)$
%for a generic nodal curve $C \s \pr^2$ of genus  0 or 1
%was found by Zariski [Zar, p. 307], and by Kaneko [Ka] for
%$d \ge 2g-1$ .
After the preprint [DeOrZa] had been
distributed the authors have received the preprint\footnote{We are thankful to
I. Shimada for sending us this preprint.} [KuShi] where a presentation
of the group $\pi_1(\pr^{2*} \se C^*)$ for a generic nodal curve
$C \s \pr^2$ of genus $g$ and degree $d \ge 2g+1$ has been computed
(see Remark 4.6 below).

It is our pleasure to thank D. Akhiezer, E. Artal, H. Flenner, V. Guba, S.
Kosarev, V. Lin, V. Sergiescu for their friendly assistance and contributions
to the paper.

\section{Big  groups and C-hyperbolicity}

\noin {\bf 1.1. Generalities on big groups}

\ss

\noin By a theorem of von Neumann, a big group is non--amenable.
The converse is not true, in general; the corresponding
examples are due to  A. Ol'shanskij, S. I. Adian and M. Gromov
(see [OSh]).
Note that the group in all these examples
is not finitely presented. For a finitely presented
group the
equivalence of bigness and non--amenability is unknown\footnote
{We are grateful to V. Sergiescu and V. Guba for this information.}.
Being non--amenable, a
big group can not be almost nilpotent or even almost
solvable. As follows from the Nielsen--Schreier Theorem, a subgroup of
finite index of a big group is big, as well
as a normal subgroup with a solvable quotient. Clearly, a group with
a big quotient is big.

\ms

We remind several classical examples of big groups. First of all, for
$g \ge 1$ the Siegel modular
group $Sp_{2g}(\gz)$ is big. In addition, it has no infinite
normal solvable subgroup (see (1.3)-(1.4) below).

\ss

Another examples are: the Artin group $B_{d, \,g}$ of the $d$--string
braids of a
genus $g$ compact Riemann surface $R_g$, and the mapping class group
Mod$_{g,\,n}$, i.e. the group of classes of isotopy of
orientation preserving diffeomorphisms of a genus $g$ Riemann surface with
$n$ punctures (see e.g. [Bi]). Namely, we have the following

\vs

\noin {\bf 1.2. Lemma.}

\noin {\it ($a$) If $g\ge2$ then $\pi_1(R_g)$ is big.

\noin ($b$) The braid group $B_{d, \,g}\,\,(d \ge 1)$ is big iff $(d,\,g)
\neq (1,\,0),\,(2,\,0),\,(3,\,0),\,(1,\,1)$.

\noin ($c$) The mapping class group Mod$_{g,\,n}$ is big iff $g \ge 1$, or
$g = 0$ and  $n \ge 4$.}

\vs

\noin {\it Proof.} ($a$)
By a theorem of Magnus [CoZi, (2.5.1)], after removing any of the standard
generators $a_1,\,b_1,\dots,\,a_g,\,b_g$ of
$\pi_1(R_g)$, the subgroup generated by the remaining ones is the free
group {\bf F}$_{2g - 1}$.

($b$)
By definition, $B_{d, \,g} =
\pi_1 (S^dR_g \setminus \De_{d, \,g})$, where $S^dR_g$ denotes the $d$-th
symmetric power of $R_g$ and $\De_{d, \,g} \subset S^dR_g$ denotes the
discriminant hypersurface consisting of the $d$-tuples of points with
coincidences. The pure braid group $P_{d,\,g} := \pi_1((R_g)^{\,d} \setminus
D_{d,\,g})$, where $D_{d,\,g} \subset (R_g)^{\,d}$ is the union of diagonal
hypersurfaces, is the normal subgroup of
$B_{d, \,g}$ of index $d$! which corresponds to the Vieta covering
$(R_g)^{\,d} \setminus D_{d,\,g} \to S^dR_g \setminus \De_{d, \,g}$.
The fibration $(R_g)^{\,d+1} \setminus D_{d+1,\,g} \to (R_g)^{\,d} \setminus
D_{d,\,g}$
with the fibre $R_g \setminus \{d\,\,\,{\rm points}\}$ yields the short
exact sequence [Bi, sect. 1.3]
$${\bf 1} \to \pi_1(R_g \setminus \{d\,\,\,{\rm points}\})
\to P_{d+1,\,g} \to P_{d,\,g} \to {\bf 1}\,.$$
For $d > 0$ the  group
$\pi_1(R_g \setminus \{d\,\,\,{\rm points}\})$ is a free group
{\bf F}$_k$ with $k = 2g + d - 1$ generators. For $d = 0$
see (a).

Hence, under
the above restrictions the pure braid group $P_{d,\,g}$, and therefore also the
braid group $B_{d, \,g}$,
contains a subgroup isomorphic to
a non-abelian free group. In the exceptional cases when $(d,\,g) = (1,\,1)$ or
$g = 0,\,1 \le d \le 3$ the same exact sequence
shows that the corresponding group $B_{d, \,g}$ is not big. This proves ($b$).

($c$) There is a natural
surjection $j\,:\,$ Mod$_{g,\,n} \to $ Mod$_g :=$ Mod$_{g,\,0}$, where the
kernel Ker$\,j$ is
the braid group $B_{n, \,g}$ if $g \ge 2$ and its
quotient by the center if $g = 1,\,n\ge 2$ or $g = 0,\, n\ge 3$
[Bi, Theorem 4.3]. Therefore, the group Mod$_{g,\,n}$ is big as soon as
the corresponding braid group $B_{n, \,g}$ is so.

For $g \ge 1$
the induced representation of Mod$_{g}$ into the first homology group of $R_g$
yields  a surjection Mod$_{g}\to$ Sp$_{2g}(\gz)$ (actually, Mod$_{1} \cong
GL(2,\, \gz)$).
This shows that Mod$_{g},\,g \ge 1$, is a big group.

For $g = 0$ we have that Mod$_{0,\,3} = B_{3, \,0}/$(center) is a finite
group, the groups Mod$_{0,\,0}$ and Mod$_{0,\,1}$ are trivial, whereas
Mod$_{0,\,2} = \gz/2\gz$ [Bi, Theorem 4.5]. This completes the proof.  \qed

\vs

\noin {\it Remark.} In fact, the Tits alternative holds in
Mod$_{g},\,g  \ge 1$ [Iv, MC] (note that for
$g \ge 2$ the latter group is not isomorphic to any arithmetic linear group
[Iv]). Furthermore,  for $g\ge 2$ any almost solvable subgroup of Mod$_{g}$ is
almost abelian [BiLuMC].

\vs

Let us make certain remarks concerning a conjectural relation of bigness
of the fundamental group of a complex space $X$ and its C-hyperbolicity.
Recall that $X$ is said to be {\it (almost) C-hyperbolic} if it has an
{\it (almost) Carath\'eodory hyperbolic} covering $Y \to X$, i.e. such that
the bounded holomorphic functions on $Y$  separate  points of $Y$
(up to finite subsets). As follows from Lin's Theorem
[Lin, Theorem B], {\it the fundamental group of an almost C-hyperbolic
algebraic variety can not be almost nilpotent.} Note that for $g \ge 1$ the
complement $\pr^{2*} \se C^*$ is C--hyperbolic, and it is almost C-hyperbolic
if $C$ is a generic rational curve of degree $d \ge 5$
[DeZa1, Thm. 1.1]. Thus, by Lin's Theorem, in all these cases
the group $\pi_1(\pr^{2*} \se C^*)$ is not almost nilpotent. Actually,
by Theorem 0.2 above this group is big. This leads us to
the following

\ms

\noin {\bf Question.} {\it Let $X$ be an almost C--hyperbolic
algebraic variety. Is then necessarily $\pi_1(X)$ a big group?}

\ms

By another theorem of Lin [Lin, Thm. B($b$)], $\pi_1(X)$
cannot be an amenable group with a non--trivial center assuming
that the universal covering space $\tilde X$ is Carath\'eodory hyperbolic.

An easy
observation is that the answer is `yes' for dim$\,X=1$. Indeed,
an algebraic curve $C$ is C--hyperbolic iff it is hyperbolic, or, in turn,
iff its normalization $C_{\rm norm}$ has a
non-abelian fundamental group. In the latter case the group
$\pi_1(C_{\rm norm})$ is big (see 1.2($a$)).
Note, however, that by a result of [LySu], any compact Riemann
surface $R$ of
genus $g \ge 2$ admits a Galois covering $\tilde R$ with a metabelian (i.e.
two-step
solvable) Galois group such that $\tilde R$ carries a non--constant bounded
holomorphic function. Modifying this result, one may even assume $\tilde R$
being Carath\'eodory hyperbolic [LinZa, Sect. 3].

\ms

More generally, we have the following fact. Its proof given below was
communicated to us by D. Akhiezer\footnote{and it is placed here with his
kind permission.}.

\vs

\noin {\bf 1.3. Theorem.} {\it Let $D \s \cz^n$ be  a bounded
symmetric domain,  and let  $\Ga \s
{\rm Aut}\,D$ be a discrete subgroup. If the Bergman volume of a
fundamental domain of $\Ga$  is finite, then $\Ga$ is a
big group and it has no infinite solvable normal subgroup.}

\vs

\noin {\it Proof.} According to a result of A. Borel and J.--L. Koszul
[Bo, Kos], a homogeneous domain $D$ is symmetric iff
the identity component $G$ of the automorphism group Aut$\,D$ is semisimple.
Recall that $G$ has trivial center, and therefore it is a connected linear
group [He, Ch. VIII.6]. Being semisimple
$G$ is not solvable. Moreover, since $G$ is connected, it is not small.  We
have  $D \cong G/K$, where $K \s G$ is a maximal compact subgroup
[ibid, Ch. VIII. 7]. The automorphism group Aut$\,D$ has finitely many
connected components, i.e. [Aut$\,D : G] < \infty$ (indeed, being a compact Lie
group the stabilizer
Stab$_z \s$ Aut$\,D$ of a point $z \in D$ has
a finite number of connected components, which is the same as those of
Aut$\,D$, because the quotient $D \simeq$ Aut$\,D/$Stab$_z$ is connected).
Hence,  $\Ga \cap G$ has finite index in $\Ga$, and the Bergman volume
of $(\Ga \cap G) \se D$ is finite, too. Therefore, the invariant volume
Vol$\,((\Ga\cap G) \se G)$ is finite, and so $\Ga\cap G$ is a lattice of $G$.

Fix a faithful linear representation $G \hookrightarrow GL(n,\,\cz)$. Let
$G_{\cz}$ be the Zariski closure of $G$ in $GL(n,\,\cz)$. By Borel's Density
Theorem (see e.g. [Ra, 5.16]), the conditions
"$G$ is semisimple and Vol$\,((\Ga\cap G) \se G) < \infty$" imply that the
subgroup
$\Ga\cap G$ is Zariski dense in $G_{\cz}$. Hence, if $\Ga$ is almost
solvable,
$G_{\cz}$ should be also almost solvable, which is not the case. By the
Tits alternative, $\Ga$ must be big.

The last assertion follows from a theorem of V. Gorbatsevich [GoShVi,
Proposition 3.7]. According to this theorem, the lattice $\Ga\cap G$ in a
connected Lie group $G$
possesses no infinite solvable normal subgroup iff $G$ is reductive and
its
semisimple part has a finite center. It is easily seen that in our case
both conditions are fulfilled. \qed

\vs

\noin {\bf 1.4.} {\it Remark.}  In fact, it would be enough in the
above theorem that $\Ga$ was a Zariski dense subgroup
of a semisimple linear algebraic group $G$ with a finite center, which acts
holomorphically in $D$. This  may be illustrated by the following
example 1.5($a$).

\ms

\noin {\bf 1.5.} {\bf Examples.}

\ss

\noin ($a$) Let $D = {\cal Z}_g$ be the Siegel upper half--plane
and $\Ga = $Sp$_{2g}(\gz),\,\,G =$ Sp$_{2g}(\rz),\,g \ge 1$, are resp.
the Siegel modular
group and the simplectic group. Then $\Ga \se D$ is a coarse moduli space
of principally polarized abelian varieties of dimension $g$, which is
a quasiprojective variety. Here $\Ga$
is Zariski dense in $G$. Actually, by a theorem of A. Borel and
Harish--Chandra [BoHC, Thm. 7.8], the arithmetic subgroup $G_{\gz}$ of a
semisimple real algebraic group $G_{\rz}$ defined over $\bf Q$ is a lattice
in $G_{\rz}$,
and so by Borel's Density Theorem, it is Zariski dense in $G_{\cz}$.
(By the way, these arguments show that $\Ga = $Sp$_{2g}(\gz)$ is a big group
without infinite normal solvable subgroups.)

\ms

\noin ($b$) Let, furthermore, $D  := T_{g,\,n}  \s\s \cz^{3g-3+n}$ be
the Teichm\"uller space  of the $n$--punctured genus $g$ marked
Riemann surfaces under the Bers realization, where $2-2g-n <0$. By
Royden's Theorem, $\Ga:=$ Aut$\,D$ is the Teichm\"uller modular
group, which coincides with
the mapping class group Mod$_{g,\,n}$. The quotient $\Ga \se D$ is
a coarse moduli space ${\cal M}_{g,\,n}$ of genus $g$ $\,\,n$--punctured
Riemann surfaces, which is a quasiprojective variety. By Lemma 1.2 above,
except the case when $(g,\,n) = (0,\,3)$ the group Mod$_{g,\,n}$ is big.

\ms

\noin (c) (see e.g. [Sh1, 2]).
A smooth projective surface
$S$ is called {\it a Kodaira surface} if there is a smooth fibration
$\pi\,:\,S \to B$ over a curve $B$, where both $B$ and a generic fibre $F$
of $\pi$ are of genus $\ge 2$ (usually $\pi$ is supposed being a non-trivial
deformation of $F$, but we don't need this assumption here).
It is well known that
the universal covering $\tilde S$ of $S$ can be realized as a bounded
pseudo--convex Bergman domain in $\cz^2$. Thus, the projective surface
$S = \Ga \se D$ is C--hyperbolic; clearly, $\Ga \simeq \pi_1(S)$ is a big
group.
More generally, the same is true when both $B$ and $F$ are quasiprojective
hyperbolic curves.

\vs

Next we pass to the simplest examples of reducible plane projective curves
with a big fundamental group of the complement.

\ms

\noin {\bf 1.6.} {\bf Examples.}

\ss

\noin ($a$) Let $C \subset \pr^2$ be a finite line
arrangement. If these lines are in general position, then by the
Deligne--Fulton Theorem, $\pi_1(\pr^2 \setminus C)$ is abelian. Otherwise,
this group is big. Indeed, let $C$ has a point $A$ of
multiplicity at least $3$. The union $L$ of lines in $C$ passing through $A$
contains at least three members of the associated linear pencil. The linear
projection $\pr^2 \setminus C \to \pr^1 \setminus \{3\,\,{\rm points}\}$ with
center at $A$ yields an epimorphism of the fundamental groups. Thus,
$\pi_1(\pr^2 \setminus C)$ dominates the free group {\bf F}$_2 = \pi_1(\pr^1
\setminus \{3\,\,{\rm points}\})$, and therefore, it is big.

In particular, if $C$ is an arrangement of six lines with four triple points,
then $\pi_1(\pr^2 \setminus C)$ is a finite index subgroup of
the mapping class group Mod$_{0,\,5}$ (see [DeZa1, 6.1($a$)];
cf. also 1.5($b$) above).

\ms

\noin ($b$) Consider further a configuration $C \subset \pr^2$ of a plane conic
together with two of its tangent lines (cf. [DeZa1, 6.1($b$)]).
The Zariski--van Kampen method yields a presentation
$$G := \pi_1(\pr^2 \setminus C) = \langle\,a,\,b\,|\,abab = baba \,\rangle\,.$$
The following proof of the bigness of $G$ was communicated to us by V.
Lin\footnote{We are grateful to
V. Lin for a kind permission to place here this proof.}.

Remind that the Coxeter group ${\rm {\bf B}}_k$ is the group generated by
the orthogonal reflections
in $\rz^k$ with respect to the coordinate planes and the diagonals $x_i - x_j
= 0,\,i,j=1,\dots,k$. The corresponding Artin--Brieskorn braid group is
the fundamental group $\pi_1(G_k({\rm {\bf B}}_k))$ of the domain $$G_k({\rm
{\bf B}}_k) := \{z = (z_1,\dots,z_k) \in \cz^k\,|\,d_k(z)\cdot z_k \neq 0\}\,,
$$ where $d_k(z)$ is the discriminant of the universal polynomial $p_k(t) =
p_k(t,\,z):= t^k + z_1t^{k-1} +\dots +z_k$
of degree $k$. Put $G_k := \{z \in \cz^k\,|\,d_k(z) \neq 0\}$, and let $E^1_k
\to G_k$ be the standard $k$--sheeted covering over $G_k$, where
$$E^1_k := \{(z,\,\la) = (z_1,\dots,z_k,\,\la) \in \cz^{k+1}\,|\, p_k(\la,\,z)
= 0\}\,.$$ Define a mapping $\,\ph\,:\,E^1_{k+1} \to G_k({\rm {\bf B}}_k)
\times \,\cz$ as follows:
$$\ph (z_1,\dots,z_{k+1},\,\la) = (q_k,\,\la) = (\xi_1,\dots,\xi_k,\,\la)\,,$$
where
$$q_k = q_k(t,\,\xi) = t^k + \xi_1t^{k-1} +\dots+\xi_k:= p_{k+1}(t+\la,\,z)/t
\in G_k({\rm {\bf B}}_k) \,.$$
Note that $t\,|\,p_{k+1}(t+\la,\,z)$, because
$p_{k+1}(\la,\,z) \equiv 0$ for $(z,\,\la) \in E^1_{k+1}$. Since
$p_{k+1}(t+\la)$ is a polynomial with simple roots, the same is true for
$q_k(t)$. Moreover,
$q_k(0) \neq 0$; thus, indeed, $q_k \in G_k({\rm {\bf B}}_k)$. It is easily
seen
that $\ph$ is a biregular isomorphism. Hence, the isomorphism
$$\pi_1(G_k({\rm {\bf B}}_k)) \cong \pi_1(E^1_{k+1}) \hookrightarrow
\pi_1(G_{k+1})$$
represents the Artin--Brieskorn braid group $\pi_1(G_k({\rm {\bf B}}_k))$ as a
subgroup of finite
index (equal to $k+1$) of the standard
Artin braid group\footnote{From now on we denote $B_m = \pi_1(G_m)$ the
standard Artin braid group with $n$ strings; don't confuse with the Coxeter
group ${\rm {\bf B}}_k$.} $B_{k+1} := \pi_1(G_{k+1})$. Therefore, the former
group is
big as soon as the latter one is so. Both of them are big starting with $k = 2$
(for the Artin group $B_{k+1}$ this can be checked in the same way as
it was done in the proof of Lemma 1.2 for the braid groups $B_{k,\,g}$). It
remains to note that $\pr^2 \setminus C \cong G_2({\rm {\bf B}}_2)$, and
therefore $G = \pi_1(\pr^2 \setminus C)$ is isomorphic to the braid group
$\pi_1(G_2({\rm {\bf B}}_2))$ which is big.

\section{Nodal approximation  of immersed curves}

Due to Theorem 2.1 below, Theorem 0.2 can be reduced
to the case where $C$ is a generic nodal Pl\"ucker curve. We also believe that
Theorem 2.1 has an independent interest.

We use the following notation and terminology.
Let $\pr^N,\,N = N(d) = {d+2 \choose 2} - 1,\,d \ge 1$, be the
Hilbert scheme of degree $d$ plane curves. Denote $Imm_{d,\,g}$ the locus of
points of $\pr^N$ which correspond to reduced irreducible immersed curves of
geometric genus $g,\,\, 0 \le g \le {d-1 \choose 2}$, and by $Nod_{d,\,g}$
resp. $PlNod_{d,\,g}$ the subset of points of $Imm_{d,\,g}$
which correspond to the nodal resp. to the Pl\"ucker nodal curves. Remind that
an irreducible curve $C \s \pr^2$ is called {\it Pl\"ucker} if the only
singular points of $C$ and the dual curve $C^*$ are ordinary nodes and cusps.
Let $Pl\ddot uNod_{d,\,g} \s PlNod_{d,\,g}$ be the subset of curves
which have no flex at a node.

Denote $\S_d \to \pr^N$ the universal family of curves over the
Hilbert scheme $\pr^N$, and let $\S_{d,\,g} \to Imm_{d,\,g}$ be its
restriction to $Imm_{d,\,g}$.  By {\it a family of curves} we  mean a proper
morphism
$\ph\,:X \to Y$ of relative dimension one of quasiprojective varieties. If
$X,\,Y$ are smooth and $\ph$ is a submersion, then the family $\ph$ is called
{\it smooth}. We say that $\ph$ admits {\it a simultaneous normalization} if
$Y$ is smooth and there exists a smooth family
of curves $\ph'\,:\,X' \to Y$ and a morphism $f\,:\,X' \to X$ commuting with
the projections onto $Y$ such that for every point $y \in Y$ the restriction
$f\,|\,X'_y\,:\,X'_y \to X_y$ onto the fibre over $y$ is a normalization map.

\vs

\noin {\bf 2.1. Theorem.} {\it

\noin a) $Imm_{d,\,g} \s \pr^N$ is a smooth locally closed subvariety of
pure dimension $3d + g - 1$.

\ss

\noin b) The universal family of curves $\S_{d,\,g} \to Imm_{d,\,g}$ admits
a simultaneous normalization $f\,:\, \M_{d,\,g} \to \S_{d,\,g}$.

\ss

\noin c) $Nod_{d,\,g}$ and, for $n \ge 2g-1,\,\,\,PlNod_{d,\,g}$ are
Zariski open subsets of $Imm_{d,\,g}$.}

\vs

\noin {\it Remark.} The first statement of (c) and the dimension count in ($a$)
can be found
in [Ha, Sect. 2], while the proofs are quite different. Note that, by Harris
[Ha], the variety $Imm_{d,\,g}$ is irreducible; it is non--empty for any
$(d,\,g)$ with $0 \le g \le {d - 1 \choose 2}$ [Se, sect.11, p.347; Ha; O].

\ms

In this section we prove ($a$), ($b$) and the first part of  ($c$) of Theorem
2.1; the proof of ($c$) is completed in sect. 3.
First we study $Imm_{d,\,g}$ locally,
in a neighborhood of a given curve $C \in Imm_{d,\,g}$.
This needs certain preparation, including a portion of plane curve
singularities.

\vs

\noin {\bf 2.2. The Gorenstein--Rosenlicht invariant, the boundary braid and
its algebraic length}

\ms

\noin Recall that the Gorenstein--Rosenlicht invariant
$\d_P$ of a singular
analytic plane curve germ $(A,\,P)$ can be expressed as $\d_P = {1 \over 2}(\mu
+ r - 1)$, where $\mu$ is the Milnor number and $r$ is the number of
local branches of $A$ at $P$ [Mi, sect. 10]. For a reduced curve $F$ on
a smooth
surface $W$ we set $\d (F) = \sum_{P \in {\rm Sing} F} \d_P$.
If $F$ is a complete irreducible curve, then by the genus formula and the
adjunction formula [BPVV, II.11] we have
\be
\pi_a (F) = g(F) + \d (F) = 1/2\, F (K_W + F) + 1\,,
\ee
where $\pi_a$ resp. $g$ denotes arithmetic resp. geometric genus,
$K_W$ is the canonical divisor of
$W$, and where for a non--compact surface $W$ we put $F K_W = {\rm
deg}\,(K_W\,|\,F)$.

\ss

Let $U \s \cz$ be the unit disc, $\Si = U \times \cz \s \cz^2$ be the solid
cylinder $\Si = \{(u,\,v) \in \cz^2\,|\,|u| < 1\}$, and $p\,:\,\cz^2 \to \cz$
be the first projection. Let $A \s \Si$ be an
analytic curve extendible transversally through the boundary $\p \Si$, so
that {\it the link} $\p A = {\bar A} \cap \p \Si$ is smooth.
Suppose also that the projection $p\,:\,A \to U$ is proper, i.e. it is a
(ramified)
covering over the unit disc $U$ of degree, say, $m$. The link $\p A$
carries a (closed) braid with $m$ strings $b_A \in B_m$ defined uniquely up
to conjugation, where $B_m$ is the Artin braid group (see (1.6($b$) above)
\footnote{
To define the braid $b_A$, cut the cylinder
$\p U = S^1 \times \cz$ along
its generator $1 \times \cz$ and then identify it with
$[0,\,1] \times \rz^2 \s \rz^3$. Fix a numbering of the points of
the fibre of $\p A$ over $1 \in \p U$. Passing once along the circle
$S^1 = \p U$ counterclockwise, we obtain the braid  $b_A $.}.

Let
$\sigma_1,\dots, \sigma_{m-1}$ be the standard generators of $B_m$.
For a braid $b = \sigma_{i_1}^{\alpha_1}\dots
\sigma_{i_n}^{\alpha_n} \in B_m$ its {\it algebraic length}
is defined as $l($b$) := \sum_{k=1}^n \alpha_k$.

\vs

\noin {\bf 2.3. Lemma.} {\it Let $A \s \Si$ as above be a nodal curve with $\d$
nodes. Suppose that all the ramification points of the covering $p\,:\,A \to U$
are simple (i.e. with ramification indices $2$) and no two of them are at the
same
fibre. If the branching divisor $D \s U$ consists of $\d + \tau$ points,
then $$l(b_A) = 2\d + \tau\,.$$}
{\it Proof.} Choose small disjoint discs $\omega_i$ in $U$,
$i=1,\dots,\d + \tau$, centered at the points of
$D$. Fix a point at the boundary of the disc $\omega_i$ and join it by
a path $\gamma_i$ with the point $1 \in \p U$, where  $\gamma_i,\,
i=1,\dots,\d + \tau$, are disjoint. The complement $U \se \bigcup_{i=1}^{\d +
\tau} ({\bar \omega}_i \cup \gamma_i)$ being simply connected, the braid $b_A$
is the product of the local braids $b_{A_i}$ which correspond to
the curves $A_i := A \cap p^{-1}(\omega_i)$. It is easily seen that the local
braid which corresponds to a node of $A$ is conjugate in the braid group $B_m$
with the square of a generator, and those at an irreducible ramification point
is conjugate with a generator. Now  the lemma easily follows. \qed

\vs

With each plane curve singularity $(A,\,{\bar 0}) \s (\cz^2,\,{\bar 0})$
we associate its {\it braid} $b_{A,\,{\bar 0}}$ defined as follows. Fix
a generic linear projection $p\,:\,(\cz^2,\,{\bar 0}) \to (\cz,\,{\bar 0})$,
so that the direction
of $p$ is different from the tangent directions of the branches of $A$ at
$\bar 0$,
and proceed in the same way as above.

\vs

\noin {\bf 2.4. Lemma.} {\it Suppose that
$(A,\,{\bar 0}) \s (\cz^2,\,{\bar 0})$ is an
immersed singularity (i.e. a singular point of a reduced curve having only
smooth local branches $A_1,\dots, A_r$) with
the Gorenstein--Rosenlicht invariant $\d = \d(A, \,{\bar 0})$.
Then

\noin a) $\,\,\,\,\,\,\,\,\,\,\,\,\,\,\,\,\d = {1 \over 2} \,
l(b_{A,\,{\bar 0}})\,.$

\ss

\noin b) Let $\tilde A$ be a small nodal deformation
of $A$ defined in a fixed small ball $B_{\epsilon}$ centered at the origin.
Denote by $r$ resp. $\tilde r$ the number of irreducible components of
$A$ resp. of $\tilde A$ in $B_{\epsilon}$.
Then $\d ({\tilde A}) \le \d (A)$, and $\d ({\tilde A}) = \d (A)$ iff $r =
{\tilde r}$. In the latter case the irreducible components ${\tilde A}_1,\dots,
{\tilde A}_r$ of ${\tilde A}$ in $B_{\epsilon}$ approximate the corresponding
irreducible components $A_1,\dots, A_r$ of $A \cap B_{\epsilon}$. }

\vs

\noin {\it Proof.} ($a$) We have $\d = \sum_{1\le k < l \le r}
(A_k\,\cdot\,A_l)_{\bar 0}$ [Mil, (10.20)]. Let $A_i' \s B_{\epsilon}$ be
a small
generic deformation of the branch $A_i,\,\,i=1,\dots,r$. Set
$A' = \bigcup_{i=1}^r A_i'$. Then $A'$ is a nodal curve with
$$\d =  \sum_{1\le k < l \le r} A'_k \cdot A'_l =
\sum_{1\le k < l \le r} (A_k\,\cdot\,A_l)_{\bar 0}$$ nodes, and
clearly, $b_{A,\,{\bar 0}} = b_{A',\,{\bar 0}}$. Since  for all $i=1,\dots,r$
the generic linear projection $p\,:\,A_i \to U_{\epsilon'}$ is
non--ramified, the same is true for the branches
$A'_i,\,i=1,\dots,r$. Thus,  $p\,:\,A' \to U_{\epsilon'}$ is ramified only
at $\d$ nodes, and therefore, in the notation of Lemma 2.3, $\tau = \tau (A')
=0$. By this lemma, we have $\d = 1/2\, l(b_{A',\,{\bar 0}})
= 1/2\, l(b_{A,\,{\bar 0}})$. This proves ($a$).

\ss

\noin ($b$) Once again here
$b_{A,\,{\bar 0}} = b_{{\tilde A},\,{\bar 0}}$.
Due to ($a$) and to Lemma 2.3, we have
$$2\d(A) = l(b_{A,\,{\bar 0}}) =
l(b_{{\tilde A},\,{\bar 0}}) = 2\d ({\tilde A}) + \tau ({\tilde A})\,,$$
and the inequality of ($b$) follows. The equality holds iff $ \tau ({\tilde A})
=0$, which means that the projection $p \,:\,{\tilde A} \to U_{{\tilde
\epsilon}}$ is ramified only at nodes of ${\tilde A}$. Therefore, for
any irreducible component ${\tilde A}_i$ of ${\tilde A} \cap B_{\epsilon}$
the composition of the normalization map
$({\tilde A}_i)_{\rm norm} \to {\tilde A}_i$ with the projection
$p \,:\,{\tilde A}_i \to U_{{\tilde \epsilon}}$ is non--ramified and hence,
one--sheeted. It follows that both of these mappings are biholomorphic,
so that the irreducible components  ${\tilde A}_i$ of ${\tilde A} \cap
B_{\epsilon}$ are smooth. The degree of the
branched covering $p\,:\,{\tilde A} \to U_{{\tilde \epsilon}}$ being equal to
$r,\,\,{\tilde A} \cap B_{\epsilon}$ consists of $r$ smooth irreducible
components close to those of $A$. \qed

\vs

Let $X$ be a smooth projective surface, $C \s X$ be an irreducible immersed
curve with a normalization $\varphi_0\,:\,M_0 \cong C_{\rm norm} \to C$. By
[No, (1.8)-(1.12)], there exists a smooth open complex surface $V$ which
contains $M_0$ as a closed subvariety, and a holomorphic immersion
$\varphi\,:\, V \to X$ that extends $\varphi_0$; it is called {\it a tubular
neighborhood of} $\varphi_0$. To obtain $V$ one simply normalizes
$C$ together with a tubular neighborhood of $C$ in $X$.

\vs

\noin {\bf 2.5. Lemma.}  {\it Let $C \s X$, $M_0$ and $V$ be as above,
and let $N \to M_0$ be the normal bundle of
$M_0$ in $V$. Then
\be {\rm deg}\, N = M_0^2 = C^2 - 2\d (C)\,.\ee
If $X = \pr^2$ and $C \in Imm_{d,\,g}$, then
\be{\rm deg}\, N = 3d +2(g-1)\,.\ee}

\noin {\it Proof.} By the adjunction formula, we have
\be 2g-2 = C^2 + CK_X - 2 \d (C) = M_0^2 + M_0 K_V \,.\ee
Since $K_V = \varphi^* K_X$, by the projection formula we have $M_0 K_V =
CK_X$, and
so (2) follows. (3) is a corollary of (2) and the genus formula (1). \qed

\vs

\noin {\bf 2.6. Corollary.} {\it a) $N$ is very ample iff $C^2 - 2\d (C)
\ge 2g+1$. For $X = \pr^2$ and $C \in Imm_{d,\,g}$ this is always the case,
and furthermore, $h^1 (M_0,\,{\cal O}(N)) = 0$
and $h^0 (M_0,\,{\cal O}(N)) = 3d + g - 1$.

\ss

\noin b) For any pair of points $p_1,\,p_2 \in M_0$ the line bundle
$N_{p_1,\,p_2} = N - [p_1] - [p_2]$ on $M_0$ is spanned\footnote{i.e. the
linear system $|N_{p_1,\,p_2}|$ has no base point.} if $C^2 - 2\d (C)
\ge 2g+2$. In particular, this is so if $X = \pr^2$ and $C \in Imm_{d,\,g}$,
where $d \ge 2$.}

\vs

\noin {\it Proof.} The first statement of (a) and (b) follow from Lemma 2.5
by the well known criteria of ampleness or spannedness
of a line bundle over a curve (see e.g. [Hart,IV.3.2] or [Na, 5.1.12]). By
the Kodaira Vanishing Theorem, we obtain that $h^1 (M_0,\,{\cal O}(N)) = 0$,
and hence, by the Riemann--Roch Formula,
we have $h^0 (M_0,\,{\cal O}(N)) = {\rm deg}\, N + 1 - g = 3d + g - 1$. \qed

\vs

The Kodaira Theorem on embedded deformations [Ko] implies such a

\vs

\noin {\bf 2.7. Corollary.} {\it There exists a maximal smooth family
$\pi_{loc}\,:\,{\cal M}_{loc} \to T_{loc}$  of embedded deformations
of the curve $M_0 \cong
\pi^{-1} (t_0)$ in $V$ over a smooth base $T_{loc}$ such that the
Kodaira--Spencer
map $T_{s_0}T_{loc} \to H^0 (M_0,\,{\cal O} (N))$ is an isomorphism. In
particular, if $X = \pr^2$ and $C \in Imm_{d,\,g}$, then\footnote{cf. [GH,
sect.2.4; Ha].} ${\rm dim}\,T_{loc} = 3d + g - 1$.}

\vs

\noin {\bf 2.8. Definition.} We say that a curve $C \in Imm_{d,\,g}$ is {\it
strongly approximated} by curves $C' \s Imm_{d,\,g}$ if $C'$ approximate $C$
in the Hausdorff topology, and for any
singular point $P$ of $C$ of multiplicity $r(C,\,P)$ and for a fixed small
neighborhood $B_{\epsilon, \,P}$ of $P$, the number $r(C',\,P)$ of irreducible
components of $C' \cap B_{\epsilon, \,P}$ is equal to $r(C,\,P)$, and the
irreducible components of $C' \cap B_{\epsilon, \,P}$ approximate those of $C
\cap B_{\epsilon, \,P}$. Or, which is equivalent, if for a
given tubular neighborhood $\ph\,:\,V \to \pr^2$ of a normalization
$\ph_0\,:\,M_0 \to C$, the curves $C'$ have  normalizations
$\ph\,|\,M'\,:\,M' \to C'$, where $M' \s V$
are obtained from $M_0$ by a small deformation.

\vs

We use below the following simple observation: a curve $C \in Nod_{d,\,g}$ is
Pl\"ucker iff $C$ has only ordinary flexes,
no multitangent line, i.e. a line tangent to $C$ in at least three points,
and no bitangent line which is an inflexional tangent. One says that
a curve $C \s \pr^2$ {\it has
only ordinary singularities} iff all the local branches of $C'$ at any of its
singular point are smooth and pairwise transversal. Denote by $Ord_{d,\,g}$ the
set of all such curves of degree $d$ and genus $g$; clearly, $Ord_{d,\,g} \s
Imm_{d,\,g}$.

\ss

The next proposition should be known at least partially; in view of the lack of
references, we give its proof.

\vs

\noin {\bf 2.9. Proposition.}  {\it The subspaces $Ord_{d,\,g},\,
Nod_{d,\,g},\,PlNod_{d,\,g}$ and
$Pl\ddot uNod_{d,\,g}$ are dense in $Imm_{d,\,g}$ in the topology of strong
approximation, and hence also in the Hausdorff topology of $\pr^N$.}

\vs

\noin {\it Proof.} Fix an arbitrary curve $C \in Imm_{d,\,g}$. We may assume
that $d \ge 3$. First we show that $C$ can be strongly approximated by curves
$C' \in Ord_{d,\,g}$.  For a curve $C' \in Imm_{d,\,g}$ denote by
$\d_1(C')$ the number of all non--ordered pairs $(A'_i,\,A'_j)$ of local
analytic branches of $C'$  which meet normally at their common center $P \in
C'$, so that $(A'_i,\,A'_j)_P = 1$. Clearly,
$C' \in Ord_{d,\,g}$ iff $\d(C') = \d_1(C')$.

Suppose that $C$ as above has a non--ordinary singular point $P$ of
multiplicity $m$. Consider the blow up $\si\,:\,X \to \pr^2$ of
$\pr^2$ at $P$, and let ${\hat C} \s X$ be the proper transform of $C$.
It is easily seen that $\d({\hat C}) = \d(C) - {m \choose 2}$. Since
${\hat C}^2 = C^2 - m^2$ we have
$${\hat C}^2 - 2\d({\hat C}) = C^2 - 2 \d(C) - m = 3d + 2(g-1) - m \ge 2g
+2\,.$$ Let $\ph \,:\,V \to X$ be a tubular neighborhood of a normalization
$\ph_0\,:\, M_0 \to {\hat C}$ of ${\hat C}$. For a pair $(A_i,\,A_j)$
of local branches of $C$ at $P$ with $(A_i,\,A_j)_P > 1$ let $\hat P \in X$
be the common center of their proper preimages ${\hat A_i},\,{\hat A_j}$
in $X$,
and let $P_i,\,P_j \in M_0$ be resp. the centers of the branches
$\ph^{-1}({\hat A_i}),\,\ph^{-1}({\hat A_j})$ of the curve $M_0 \s V$. By Lemma
2.5, for the normal bundle $N$ of $M_0$ in $V$ we have
$${\rm deg}\,(N - [P_i]) = M_0^2 - 1 = {\hat C}^2 - 2\d({\hat C}) - 1
\ge 2g + 1\,.$$ Therefore, being spanned, the line bundle $N - [P_i]$ possesses
a section which does not vanish at $P_j,\,j \neq i$. It follows that $N$ has a
section
that vanishes at $P_i$, but not at $P_j$. This yields a deformation
$M'_0$ of $M_0$ in $V$ which passes through $P_i$, but not through $P_j$. Thus,
for the curve $C' := \si\ph (M_0') \in Imm_{d,\,g}$ close enough to $C$ we have
 $\d_1(C') > \d_1(C)$. By induction on $\d_1(C')$, we get a strong
approximation $C' \in Ord_{d,\,g}$ of $C$.

Suppose further that $C \in Ord_{d,\,g}$
is not nodal, i.e. it has  a point $P$ of multiplicity $m \ge 3$.  Applying
the same procedure as above to a triple of points $P_i,\,P_j,\,P_k \in M_0$
which lie over $P$, and using the inequality $${\rm deg}\,(N - [P_i] - [P_j])
\ge 2g\,,$$ by the spannedness of the line bundle $N - [P_i] - [P_j]$,
we obtain a section of $N$ which vanishes at the points $P_i$ and $P_j$, but
not at $P_k$. This leads to a curve $C' \in Ord_{d,\,g}$ which strongly
approximates $C$
and is simpler than $C$ in the following sense: $m(C') < m(C)$, where
 $$m(C) := \sum_{P_i \in {\rm sing}\,(C)} ({\rm mult}\,(P_i) - 1)\,.$$
Induction on $m(C')$ now shows that $C$ can be strongly approximated by curves
$C' \in Nod_{d,\,g}$.

\ss

Next we show that a curve $C \in Nod_{d,\,g}$ can be strongly approximated by
curves
$C' \in Nod_{d,\,g}$ with only ordinary flexes. We proceed by induction on the
number $ofl(C')$ of ordinary flexes of $C'$. Since such a flex is a normal
intersection point of
the Hesse curve of $C'$ with a smooth local branch of $C'$, clearly,
the bounded function $ofl(C')$ is lower semi--continuous on $Nod_{d,\,g}$ with
respect to the Hausdorff topology.

Suppose that $C$ has
a non--ordinary flex at a local branch $A$ of $C$ centered at $P \in C$, so
that $(A,\,L)_P \ge 4$, where $L$ is the tangent line to $A$ at $P$. In the
notation as above, let $\si\,:\,X \to \pr^2$ be the composition of three
successive blow ups over $P$ with centers at the proper preimages of $A$.
Let ${\hat L} \s X$ be the proper transform of $L$ and $\hat P$ be the center
of the proper transform ${\hat A} \s X$ of $A$. We have
$({\hat A},\, {\hat L})_{\hat P} \ge 1$. If $P$ is a smooth point of $C$ then
${\hat C}^2 = C^2 - 3$ and $\d({\hat C}) = \d (C)$. If $P$ is a node of $C$
then ${\hat C}^2 = C^2 - 6$ and $\d({\hat C}) = \d (C) - 1$. In any case,
$${\rm deg}\,N = {\hat C}^2 - 2\d({\hat C}) \ge C^2 - 2\d(C) - 4 \ge 2g\,.$$
Therefore, the normal bundle $N$ of $M_0$ in $V$ is spanned, and hence it has a
section which does not vanish
at the point $\hat P$. The corresponding Kodaira-Spencer deformation yields a
curve $M_0'$ on $X$ close enough to $M_0$ which does not pass through $\hat P$.
It is easily seen that the projection $C' := \si\ph(M_0') \s \pr^2$ is a nodal
curve  with an ordinary flex at $P$ and such that $ofl(C') > ofl(C)$. After a
finite
number of steps we obtain a strong approximation $C' \in Nod_{d,\,g}$ of $C$
with only ordinary flexes.

Suppose  further that $C \in Nod_{d,\,g}$ has only ordinary flexes
(note that this  is an open condition). We will
find a strong approximation $C'$ of $C$ without multiple tangents.
Denote by $b(C')$ the total number of distinct intersection points with $C$ of
all the bitangent lines of $C'$. Clearly, the bounded function $b(C')$
is lower semi--continuous on $Nod_{d,\,g}$.

Let $C$ have a multitangent line $L$ which is tangent to $C$ at points
$P,\,Q,\,R \in C$ and, perhaps, at some other points.
Let $\si\,:\,X \to \pr^2$ be the composition of the blow-ups of $\pr^2$ at the
points $P,\,Q$ and $R$, and let $\hat C$ be the proper transform of $C$ at $X$.
Note that $d = $ deg$\,C = L\cdot C \ge 6$.  As above, the blow up at a smooth
point (resp. at a node) of $C$ decreases the difference $C^2 - 2\d (C)$ by $1$
(resp. by  $2$). Thus, we have
$${\hat C}^2 - 2\d ({\hat C}) \ge C^2 - 2\d (C) - 6 = 2g + 3d - 8 \,.$$
Let $\ph\,:\,V \to X$ be a tubular neighborhood of a normalization
$\ph_0\,:\,M_0 \to {\hat C}$ of $\hat C$, and let $N$ be the normal bundle of
$M_0$ in $V$.
The line bundle $N - [{\hat P}] - [{\hat Q}]$ on $M_0$ of degree
$\ge 2g + 3d - 10 > 2g + 2$ is spanned (cf. Corollary 2.6). This yields a
deformation $C' := \si \ph(M_0') \s Nod_{d,\,g}$ of $C$ such that $L$
is still tangent to $C'$ at the points $P$ and $Q$, and meets $C'$  normally
at $R$, so that $b(C') > b(C)$. Maximizing $b(C')$ we get a
strong approximation $ C' \in Nod_{d,\,g}$ of $C$ with only ordinary flexes
and without multiple tangents.

Suppose now that $C \in Nod_{d,\,g}$ has only ordinary flexes and no multiple
tangent line, which is an open condition. To find a strong Pl\"ucker
approximation
$C'$ of $C$, we will proceed by induction on the total number $inf(C')$ of
distinct intersection points of $C'$ with all of its inflexional tangent lines.
We have to ensure that no inflexional tangent line of $C'$ is a bitangent line.

Let a bitangent line $L$ of $C$ be an inflexional tangent of $C$ at a point
$P \in C$ and tangent to $C$ at a point $Q \in C$. Then
$d =$ deg$\,C = C\cdot L \ge 5$. Blowing up $\pr^2$ at $Q$
we get a surface $X = \si_Q(\pr^2)$. In the notation as above, we have
$${\rm deg}\,(N - 3[{\hat P}]) \ge 2g + 3d - 7 > 2g + 2\,.$$ Therefore,
there exists a deformation $ C' = \si \ph(M_0') \in Nod_{d,\,g}$ of $C$ such
that $L$ is still an inflexional tangent of $C$ at $P$, but it meets $C$
normally at $Q$. Thus, $inf(C') > inf(C)$. By induction, we obtain
a strong approximation $C'$ of $C$ which belongs to $PlNod_{d,\,g}$.

Suppose finally that $C \in PlNod_{d,\,g} \setminus Pl\ddot uNod_{d,\,g}$,
so that, although all the flexes of $C$ are
ordinary, one of them, say $(A,\,P)$, is located at a node of $C$ with the
second branch, say, $B$. This time we proceed by induction on
the number $sfl(C')$ of flexes of $C'$ which are smooth
points. Evidently,  $sfl(C')$ is a bounded lower
semi--continuous function  on $Nod_{d,\,g}$.

Performing two successive blow ups, the first one at $P \in C$ and the second
one at the center of the proper transform of the branch $A$, we obtain a
surface $X$. Denote by
${\hat Q}$ the center of the proper preimage $\hat B$ of the branch $B$
in $X$. We have
$${\rm deg}\,N = {\hat C}^2 - 2\d({\hat C}) = C^2 - 2\d(C) - 3
\ge 2g + 1\,,$$
so that the line bundle $N - [{\hat Q}]$ on $V$ is spanned.
Hence, we can find a
section of $N$ which vanishes at $\hat Q$ and does not vanish at $\hat P$.
This yields a small deformation $M_0'$ of $M_0$ on $V$ which passes through
$\hat Q$ but not through $\hat P$. The curve $C' := \si\ph(M_0') \s \pr^2$
is close enough to $C$, still has a node at $P$ which is not any more
a flex, while $L$ is still a tangent line of $C'$ at $P$. Note that a small
deformation of $C$ yields a small deformation of the Hesse curve $H_C$ of $C$,
so that the flexes of $C$ which are the (normal) intersection points of $C$
and $H_C$ are also perturbed a little. Thus, $C'$ has a flex at a smooth point
close to $P$. It follows that $sfl(C') > sfl(C)$. In a finite number of steps
we obtain a desired strong approximation $C' \in Pl\ddot uNod_{d,\,g}$ of $C$.
This completes the proof. \qed

\vs

The next lemma  shows  that the strong approximation of immersed
curves coincides with the usual one.

\vs

\noin {\bf 2.10. Lemma.} {\it Let $C \in Imm_{d,\,g}$, and let
$\varphi\,:\, V \to \pr^2$ be a tubular neighborhood of its normalization
$\varphi_0\,:\,M_0 \to C$.
Then any curve $C' \in Imm_{d,\,g}$ close enough to $C$ in
the Hausdorff topology of $\pr^N$ (or, which is the
same,  coefficientwise) is the image of a unique smooth curve
$M  \cong C'_{\rm norm}  \s V$ under the holomorphic mapping
$\varphi\,:\,V \to \pr^2$.}

\vs

\noin {\it Proof.} Let $P$ be a singular point of $C$, and let $B_{\epsilon,
\,P}$ be a fixed small neighborhood of $P$. Denote by $r(C,\,P)$ the
multiplicity of $C$ at $P$, and by $r(C',\,P)$ the number of irreducible
components in $B_{\epsilon, \,P}$ of a curve $C'$ close enough to $C$
(cf. Definition 2.8).
Once we show that $\,r(C,\,P) = r(C',\,P)$ for any singular point
$P$ of $C$, then the irreducible components
of $C' \cap B_{\epsilon, \,P}$ approximate those of $C \cap B_{\epsilon, \,P}$,
i.e. $C'$ is a strong approximation of $C$, and the statement follows.

Actually, it is sufficient to prove the equality  $r(C,\,P) = r(C',\,P)$ under
the additional assumption that the approximating curve $C'$ is nodal.
Indeed, by Proposition 2.9, the curve
$C'\s Imm_{d,\,g}$ can be, in turn, strongly approximated by a curve $C'' \in
Nod_{d,\,g}$. Since $C''$ approximates both $C$ and $C'$ in the Hausdorff
topology, from the
equalities $r(C'',\,P) = r(C,\,P)$ and $r(C'',\,P) = r(C',\,P)$ it follows
that $r(C',\,P) = r(C,\,P)$.

Assuming further that $C'$ is nodal, by (1)
and Lemma 2.4($b$), we obtain
$${n-1 \choose 2} - g = \d (C') = \sum_{P\in {\rm Sing}\,C'} \d (C'\cap
B_{\epsilon, \,P}) \le \sum_{P\in {\rm Sing}\,C} \d(C,\,P) =
{n-1 \choose 2} - g\,.$$ Henceforth,
$\d (C'\cap B_{\epsilon, \,P}) = \d (C, \,P)$ for all
$P\in {\rm Sing}\,C$.
Applying Lemma 2.4($b$) once again, we get that  $r(C',\,P) = r(C,\,P)$
for all $P\in {\rm Sing}\,C$, as desired. \qed

\vs

\noin {\bf 2.11. Lemma.} {\it ($a$) $Imm_{d,\,g}$ is a locally closed complex
analytic submanifold of $\pr^N$ of dimension $3d + g - 1$.

\ss

\noin ($b$) The universal family of curves $\S_{d,\,g} \to Imm_{d,\,g}$
over $Imm_{d,\,g}$ admits a complex analytic simultaneous normalization
$f = f_{d,\,g}\,:\,\M_{d,\,g} \to \S_{d,\,g}$.}

\vs

\noin {\it Proof.} Fix a curve $C \in Imm_{d,\,g}$, and consider a tubular
neighborhood $\ph\,:\,V \to \pr^2$ of a normalization $\ph_0$ of $C$. By
Corollary 2.7 and Lemma 2.10, the projection $\ph$ yields a local
analytic chart $U_C$ of dimension $3d + g - 1$ on $Imm_{d,\,g}$ centered at
$C$ which covers the whole intersection of $Imm_{d,\,g}$ with a sufficiently
small ball in $\pr^N$ around $C$. This proves ($a$).

\ss

To prove ($b$) denote by $\S_C$ the restriction of the family $\S_{d,\,g}$ onto
the chart $U_C$. Note that the same projection $\ph$ yields an
analytic simultaneous normalization $f_C\,:\,\M_C \to \S_C$ of $\S_C$. Any two
such normalizations $f_C\,:\,\M_C \to \S_C$ and $f'_C\,:\,\M'_C \to \S_C$ over
the
same chart $U_C$  which  arise  from two different tubular neighborhoods
$\ph,\,\ph'$, can be naturally biholomorphically identified via their
projections. Hence, the equivalence class of these normalizations over the
same chart $U_C$ in $Imm_{d,\,g}$ can be regarded as an equivalence
class of  charts  on a new complex manifold $\M_{d,\,g}$ of dimension $3d + g$.
Indeed, suppose that two charts $U_C$ and $U_{C'}$ on $Imm_{d,\,g}$ have a
non--empty intersection $U_{C,\,C'} := U_C \cap U_{C'}$. Consider a fibrewise
bimeromorphic mapping of smooth manifolds
$f_{C,\,C'} := f_{C'}^{-1} \circ f_C\,:\, \M_C\,|\,U_{C,\,C'} \to
\M_{C'}\,|\,U_{C,\,C'}$. It is biholomorphic at the complement of the
`multiple point locus'
$D_{C,\,C'} := f_C^{-1}($sing$\,S_{C,\,C'})$, where $S_{C,\,C'}:= S_C\,|\,
U_{C,\,C'}$, and by Riemann's extension Theorem, it has a holomorphic extension
through $D_{C,\,C'}$.  Clearly, the projection
$f_{d,\,g}\,:\,\M_{d,\,g} \to \S_{d,\,g}$ induced by the local mappings
$f_C\,:\,\M_C \to \S_C$ is a holomorphic simultaneous normalization, which
proves ($b$). \qed

\vs

Next we show that all the above subvarieties of the Hilbert scheme $\pr^N$
are algebraic. Although
the following statement holds in much bigger generality\footnote{We are
grateful to H. Flenner who introduced to us this circle of ideas.} (cf. e.g.
[BinFl, Theorem 2.2]), it will be enough for us this
restricted version which has a rather easy proof.

\vs

\noin {\bf 2.12. Lemma.} {\it Let $f\,:\,X \to Y$ be a family of curves over an
irreducible base $Y$. Then there exists a Zariski open subset $U \s Y$ such
that the restriction $f\,|\,f^{-1}(U)$ of $f$ over $U$  admits a simultaneous
normalization.}

\vs

\noin {\it Proof.} Without loss of generality  we may suppose $Y$ being smooth.
Let $\nu\,:\,X_{\rm norm} \to X$ be a normalization. Consider the induced
family of curves $f' := f \circ \nu$. Since the singular locus $S$ of
the normal variety $X_{\rm norm}$ has codimension at least $2$, its image
$f'(S) \s Y$ has codimension at least $1$. Restricting $f$ and $f'$ onto
the complement of the Zariski closure $\overline {f'(S)}$ of the constructible
subset $f'(S)$ in $Y$, we may suppose $X_{\rm norm}$ being smooth. By the
Bertini--Sard Theorem [Hart, III.10.7], $f'$ is an immersion over a Zariski
open
subset $U \s Y$. Therefore, each fibre $(f')^{-1}(y),\,y \in U$, is smooth,
and the restriction $\nu\,|\,(f')^{-1}(y)$ yields a normalization of the curve
$X_y :=f^{-1}(y)$. Thus, we have obtained the desired simultaneous
normalization of the original family $f$ over $U$. \qed

\vs

We use below the following notation. Given a family of curves $f\,:\,X \to Y$,
for any $g \ge 0$ denote by $Curv_g(f)$ the subset of points $y \in Y$ such
that the fibre $X_y$ over $y$ is a reduced irreducible curve of geometric
genus $g$. For the universal family $f_d\,:\,\S_d \to \pr^N$ of  degree $d$
curves in $\pr^2$, set $Curv_{d,\,g} = Curv_g(f_d)$.

We say that an abstract reduced irreducible curve $C$ is {\it of immersed type}
if its normalization map $\nu\,:\,C_{\rm norm} \to C$ has a nowhere vanishing
differential. Let $Imm_g(f)$ be the subset of points $y \in  Curv_g(f)$ which
correspond to the curves of immersed type, so that, in particular, $Imm_g(f_d)=
Imm_{d,\,g}$.

\vs

\noin {\bf 2.13. Corollary.} {\it ($a$) Given a family of curves $f\,:\,X \to
Y$, the base $Y$ can be represented as a disjoint union of smooth irreducible
quasi--projective subvarieties $Y_i \s Y,\,\,i=1,\dots,n = n(f)$, such that
for
each $i=1,\dots,n $ the restriction of $f$ onto $Y_i$ admits a simultaneous
normalization.

\ss

\noin ($b$) For any $g \ge 0$ the subsets $Curv_g(f) \s Y$ and $Imm_g(f) \s Y$
are constructible. In particular, $Curv_{d,\,g}$ and $Imm_{d,\,g}$ are
constructible subsets of the Hilbert scheme $\pr^N$.}

\vs

\noin {\it Proof.} ($a$) Assuming for simplicity that $Y$ is irreducible
we start with $Y_1 := U$, where $U \s Y$ is as in Lemma 2.12 above. Next we
apply Lemma 2.12 to the restriction of $f$ onto each of the irreducible
components of the regular part of the Zariski closed subvariety $Y^{(1)}:=Y
\setminus Y_1$ of $Y$. Following this way, in a finite number of
steps we obtain the desired partition of $Y$. \qed

\ss

\noin ($b$) Since $f\,|\,Y_i$ admits a simultaneous normalization, for any
$i=1,\dots,n$ the number and the geometric genera of the irreducible components
of a fibre $X_y = f^{-1}(y)$ do not depend on $y \in Y_i$. Thus, $Curv_g(f)$ is
a union of some of the $Y_i$, and hence it is constructible.

Set $X_i = f^{-1}(Y_i)$ and $f_i = f\,|\,X_i$, where $Y_i \s Curv_g(f)$ is a
stratum of the above stratification. Let
$$
\begin{picture}(800,60)
\unitlength0.2em
\thicklines
\put(88,1){$Y_i$}
\put(64,23){$X_i'$}
\put(108,23){$X_i$}
\put(83,24){$\vector(1,0){15}$}
\put(105,19){$\vector(-1,-1){12}$}
\put(72,19){$\vector(1,-1){12}$}
\put(70,10){$p_i$}
\put(89,29){$\nu_i$}
\end{picture}
$$
be a simultaneous normalization. Denote by $T_{Y_i} X_i' = $ Ker$\,dp_i$ the
relative tangent bundle of $p_i$; $p_i$ being a smooth family of curves,
$T_{Y_i} X_i'$ is a smooth line bundle on $X_i'$. Let $D_i \s X_i'$ be the
locus of points where the restriction $d\nu_i\,|\,T_{Y_i} X_i'$ vanishes.
Since $D_i$ is Zariski closed its image $p_i(D_i) \s Y_i$ is a constructible
subset of $Y_i$. Clearly, the complement $Y_i \setminus p_i(D_i)$ coincides
with $Imm_g(f) \cap Y_i$. Thus, the latter subset is constructible for all
$i= 1,\dots,n$. Hence, $Imm_g(f)$ is constructible, too. \qed

\vs

\noin {\bf 2.14.} {\it Starting the proof of Theorem 2.1.} ($a$)
directly follows from Lemma 2.11($a$) and Corollary 2.13($b$).  From
Corollary 2.13($b$) it also follows
 that the total space $\S_{d,\,g}$ of the universal
family of curves $\S_{d,\,g} \to Imm_{d,\,g}$ is a quasi--projective variety.
The holomorphic mapping $f = f_{d,\,g}\,:\,\M_{d,\,g} \to \S_{d,\,g}$ which
realizes an analytic simultaneous normalization is finite and proper
(see Lemma 2.11($b$)). Therefore, by the Grauert--Remmert Theorem [Ha, B 3.2],
$\M_{d,\,g}$ possesses a structure
of a quasi--projective variety,  so that $f$ is a finite morphism of
quasi--projective varieties. Thus, $f$  yields an algebraic simultaneous
normalization of the universal family of curves over  $Imm_{d,\,g}$. This
proves ($b$).

To prove the first part of (c) denote $T = Imm_{d,\,g},\,\,\S_T =
\S_{d,\,g},\,\,\M_T = \M_{d,\,g}$ and $\pr^2_T = \pr^2 \times T$. There is
a natural embedding $\,i\,:\,\S_T \hookrightarrow \pr^2_T$. Consider the
composition $\ph:= i \circ f\,:\,\M_T \to \pr^2_T$ and its relative square
$\ph^{(2)} := \ph^2_T \,:\,\M_T^2 \to (\pr^2_T)^2$, where $\M_T^2 := \M_T
\times_T \M_T$ and $(\pr^2_T)^2 := \pr^2_T\times_T \pr^2_T$.  Let $\D_T \s
\M_T^2$ resp. $D_T \s (\pr^2_T)^2$ be the diagonals. Clearly, $E:=
(\ph^{(2)})^{-1}(D_T) \setminus
\D_T$ is a closed subvariety of $\M_T^2$, and the restriction
$\pi^{(2)}\,|\,E\,:\,E \to T$ of the projection $\pi^{(2)}\,:\,\M_T^2 \to T$
has finite fibres. Its fibre over a point $t \in T = Imm_{d,\,g}$ corresponds
to the multiple point divisor on the normalization $M_t$ of the immersed curve
$S_t \s \pr^2$. The restriction $\ph^{(2)}\,|\,E\,:\,E \to D_T$ is a finite
morphism. The image ${\tilde E} := \ph^{(2)}(E)$ is proper over $T$. Moreover,
the fibre $\tau^{-1}(t) \s {\tilde E}$ over a point $t \in T$ under the
restriction to ${\tilde E}$  of the projection $\tau\,:\,D_T \to T$ corresponds
to the set of singular points of the curve $S_t$. Therefore, it
consists of  $\d = {d-1 \choose 2} - g$ points iff $S_t$
is a nodal curve. By Proposition 2.9, any irreducible component of $T$ contains
points which correspond to nodal curves. Thus, the finite morphism
$\tau\,:\,{\tilde E} \to T$
has degree $\d$ over every such component, and so, the complement
$Imm_{d,\,g} \setminus Nod_{d,\,g} \s T = Imm_{d,\,g}$ coincides with the
ramification divisor $R_{\tau}$ of $\tau$. Hence, $Nod_{d,\,g} \s Imm_{d,\,g}$
is, indeed, a Zariski open subvariety. \qed

\vs

\noin {\it Remark.} If $S_t$ is a nodal curve with
$\d = {d-1 \choose 2} - g$ nodes, then the fibre $p^{-1}(t)$  of the
above projection $p := \pi^{(2)}\,|\,E\,:\,E \to T$ consists of $2\d$ points.
The latter holds true if $S_t$ has only ordinary singularities. Hence, the
subset $Imm_{d,\,g} \setminus Ord_{d,\,g}$ is contained in the ramification
divisor $R_p \s T$ of $p$.

\section{Pl\"ucker conditions}

It  is known  [Au] that in general, the subset of
the rational Pl\"ucker curves is not Zariski open
in the space ${\cal R}_d$ of all the rational plane curves of  a given
degree $d$,  although it always contains a Zariski open subset of
${\cal R}_d$.  Nevertheless, we will show that $PlNod_{d,\,g}$ is a Zariski
open subset of $Imm_{d,\,g}$ for $d \ge 2g-1$, which proves Theorem 2.1(c).

\vs

\noin {\bf 3.1. Lemma.} {\it Let $C \s \pr^2$ be an irreducible nodal curve
of
degree $d$ with the normalization $M$, and let $g^2_d$ be
the linear system on $M$
of all line cuts of $C$. Then $C$ is a Pl\"ucker curve iff $g^2_d$ contains
no divisor $D$ of the form
$$(i) \,\,\,D = 4p_1 +\dots\,;\,\,\,\,\,\,or\,\,\,\,\,\,(ii)\,\,\,D = 3p_1 +
2p_2 +\dots\,;\,\,\,\,\,\,or\,\,\,\,\,\,(iii)
\,\,\,D = 2p_1 + 2p_2 + 2p_3+\dots\,,$$ where $p_i \in M$ are not necessarily
distinct.}

\vs

\noin {\it Proof.} The system $g^2_d$ contains no divisor of type (i)
iff all the flexes of $C$
are ordinary, i.e. all the singular branches of $C^*$ are ordinary cusps.
Under this condition, at most two of the local branches of $C^*$ meet at
a point
iff $g^2_d$ does not contain any divisor of type (iii).
Furthermore, two branches of $C^*$
meet at a point and one of them is singular iff $g^2_d$ contains a
divisor $D$ as in (ii). Since
$C$ being nodal has no tacnode, $C^*$ has no one, too. Therefore, $C$ is a
Pl\"ucker curve iff $g^2_d$ does not contain any divisor $D$ as in (i)--(iii).
\qed

\vs

Recall  the following notion (see e.g. [ACGH]).

\ms

\noin {\bf 3.2. Picard bundles}

\ss

Let $M$ be a smooth projective curve of genus $g$. The $d$--th symmetric
power $S^dM$ (which is a smooth manifold) might  be regarded as the space of
degree
$d$ effective divisors on $M$. Let $J_d(M) = $Pic$\,^d(M)$ be the component
of the Picard group Pic$\,(M)$ which parametrizes the degree $d$ line
bundles on
$M$, and let $\phi_d \,:\,S^dM \to J_d(M)$ be the morphism sending a
degree $d$ effective divisor on $M$ into its linear equivalence class.
Chosing a base point $p_0 \in M$ we may identify $ J_d(M)$ with the
Jacobian variety
$J_0(M)$ and $\phi_d$ with the $d$-th Abel--Jacobi mapping.
By a theorem of Mattuck [Ma] (see also [ACGH, Ch.IV]) for $d \ge 2g -1$
the morphism $\phi_d$ is a submersion and moreover, it defines a projective
bundle (i.e. a projectivization of an algebraic vector bundle) with the
standard fibre $\pr^{d-g}$. This bundle is called {\it the $d$-th Picard
bundle of $M$}.

Given a smooth family $\pi\,:\,\M \to T$ of complete genus $g$
curves and given $d \ge 2g -1$,
there is the associated Picard bundle $\Phi_d \,:\,S^d\M \to \J_d(\M)$
of relative smooth schemes over $T$. Consider also the associated
grassmanian bundle
$Grass_{2,\,d-g}(\M) \to \J_d(\M)$ which parametrizes the two--dimensional
linear series $g^2_d$ of degree $d$ on the fibres $M_t =
\pi^{-1}(t),\,t\in T$.

Let $\pi\,:\,\M \to T,\,\,T = Imm_{d,\,g}$, be the family constructed in 2.14
above. Then for each $t \in T$ there is the linear series
$g^2_d = g^2_d(t)$ on $M_t$ of the line
cuts of the plane curve $C_t = f(M_t) \s \pr^2$. This defines a regular
section $\sigma\,:\,T \to Grass_{2,\,d-g}(\M)$.

\vs

\noin {\bf 3.3.} {\it Finishing up the proof of Theorem 2.1(c).} Let
$\pi\,:\,\M := \M_T \to T$ be the family as in 2.14, and let
$\Phi_d \,:\,S^d\M \to \J_d(\M)$
be the associated Picard bundle. Denote $\D^{(i)}$ resp.
 $\D^{(ii)},\,\,\D^{(iii)}$ the subvariety of $S^d\M$ which consists of the
degree $d$ effective divisors on the fibres $M_t$ of $\pi$ of the form
(i) resp. (ii), (iii) of Lemma 3.1. Set $\D =  \D^{(i)} \cup \D^{(ii)}
\cup \D^{(iii)}$. Note that $\D$ is a closed subvariety of $S^d\M$ of
codim$\,_{S^d\M}\D \ge 3$ (and moreover,
codim$\,_{S^dM_t}\D_t \ge 3$ for each $t \in T$). Indeed, to be in $\D_t$
a divisor on $M_t$ must satisfy a system of three independent equations.

Let $\Z \s Grass_{2,\,d-g}(\M) \times S^d\M$ be the incidence relation.
Its fibre $\Z_t$ over a point $t \in T$ consists of all pairs $(L,\,v)$, where
$L$ is a two--plane in $\pr^{d-g}_j := \phi_d^{-1}(j),\,\,j \in J_d(M_t)$, and
$v \in \pr^{d-g}_j$ is a point of $L$. Let $pr_1\,:\,\Z \to
Grass_{2,\,d-g}(\M),\,\,pr_2\,:\,\Z \to S^d\M$ be the canonical projections,
and let $\si\,:\,
T \to Grass_{2,\,d-g}(\M)$ be the regular section as in (3.2) above.
Put
$\Z_{\D} := pr_2^{-1}(\D) \s \Z$, ${\hat \D} := pr_1(\Z_{\D}) \s
Grass_{2,\,d-g}(\M)$ and $T' := \si^{-1}({\hat \D}) \s T$.
Since the projection
$pr_1$ is proper, ${\hat \D} \s Grass_{2,\,d-g}(\M)$, and therefore also
$T' \s T$ are closed subvarieties of the corresponding varieties. Clearly,
$t \in T'$ iff
the linear series $g^2_d(t) = \si(t)$ on $M_t$ contains a divisor from $\D_t$.

Recall that $Nod_{d,\,g} = T \setminus R_{\tau}$, where $R_{\tau}$ is the
ramification divisor as in (2.14). By Lemma 3.1, we have that
$PlNod_{d,\,g} = T \se (R_{\tau} \cup T')$. By Proposition 2.9, any irreducible
component $I$ of $T = Imm_{d,\,g}$ contains a Pl\"ucker curve. Thus,
$T' \cap I$ is a proper subvariety of $I$; in particular,
codim$\,_{T}T' \ge 1$. Hence, $PlNod_{d,\,g}$ is, indeed, a Zariski open
subset
of $T = Imm_{d,\,g}$. This completes the proof of Theorem 2.1. \qed

\vs

Theorem 2.1 implies

\vs

\noin {\bf 3.4. Corollary.} {\it Any irreducible plane curve $C^*$ of genus
$g$
and degree $n = 2(d + g - 1)$, where $d \ge 2g-1$, whose dual $C$ is an
immersed curve, is a specialization
of generic maximal cuspidal  Pl\"ucker curves $C'^*$ of the same degree and
genus\footnote{i.e. $C'^*$ has the maximal  number of cusps allowed by
Pl\"ucker's formulas.}. Hence, there is an epimorphism
$\pi_1(\pr^{2*} \se C^*) \to \pi_1(\pr^{2*} \se C'^*)$.
In particular, the former group is big (resp. non--amenable, non--almost
solvable, non--almost nilpotent) if the latter one is so.}

\vs

\noin {\it Proof.} By the class formula [Na, 1.5.4], the dual of an
irreducible immersed plane
curve of degree $d$ and genus $g$ has degree $d^* = 2(g + d - 1)$. By
Theorem 2.1($a$) and ($b$),
there is the diagram
$$
\begin{picture}(800,60)
\unitlength0.2em
\thicklines
\put(88,23){${\cal M}_T$}
\put(64,1){$\pr^2_T$}
\put(108,1){$\pr^{2^*}_T$}
\put(86,2){$\longleftrightarrow$}
\put(87,18){$\vector(-1,-1){12}$}
\put(93,18){$\vector(1,-1){12}$}
\put(70,14){$f$}
\put(110,14){$f^*$}
\end{picture}
$$
where the morphism $f^*$ yields a simultaneous normalization of the dual
family, so that for each $t \in T = Imm_{d,\,g}$
the image $f^*(M_t) =S_t^*$ is the dual curve of the curve $S_t = f(M_t)$
(see e.g. [Na, 1.5.1]).
By ($c$), the subset $PlNod_{d,\,g}\s T$  is  Zariski open. The dual $S_t^*$,
where $t \in PlNod_{d,\,g}$, is a
maximal cuspidal curve of degree $d$ and genus $g$. Vice versa, any such curve
is the dual $S_t^*$ of a nodal Pl\"ucker curve $S_t,\,t \in PlNod_{d,\,g}$.
This yields the first assertion. The second one follows from a well known
theorem of Zariski (see [Zar, p.131, Thm.5] or [Di, 4.3.2]). As for the third
one, see (1.1) above. \qed

\section{Proof of Theorem 0.2}

\vs

The following lemma is a particular case of the Varchenko Equisingularity
Theorem [Va, Theorem 5.3].

\vs

\noin {\bf 4.1. Lemma.} {\it Let $p\,:\,E \to B$ be a surjective morphism,
where $E$, $B$ are smooth connected quasi--projective varieties. Then
there exist a proper subvariety  $A \s B$   such that the restriction $p\,|\,(E
\se H)$,  where
$H = p^{-1}(A)$,  determines
a smooth locally trivial fibre bundle $p \,:\,E \se H \to B \se A$.}

\vs

 Let $\De$ be a hypersurface in a complex manifold $E$, $e \in$
reg$\,\De$ be  a smooth point  of $\De$, and $\o$ be  a small
disc in $E$ centered at $e$ and transversal to $\De$.
By {\it a vanishing loop} of $\De$ at  $e$
we mean a loop $\d$ in $E \se \De$ consisting of a path $\alpha$ which
joins a base point $e_0 \in E \se \De$ with a point $e' \in \o \se
\De$  and a loop $\beta$ in $\o \se \De$ with the base point $e'$
(i.e. $e$ is in the interior of $\beta$ in $\omega$).

The next simple lemma is well known; for the sake of completeness
we give its proof.

\vs

\noin {\bf 4.2. Lemma.} {\it Let, as before, $\De$ be a hypersurface in a
complex manifold $E$, and let $\g_0,\,\g_1\,:\,S^1 \to E \se \De$ be two loops
with
the base point $e_0 \in E \se \De$ joined in $E$ by a smooth homotopy $\g\,:\,
S^1 \times [0,\,1] \to E$ transversal to $\De$, such that the image $S =$
Im$\,\g$ meets $\De$ at  the points $e_1,\dots,e_k \in {\rm reg}\,\De$. Then
$\g_0$ is homotopic in $E \se \De$ to a product
$\g_1\d_{i_1}\dots \d_{i_k}$, where $(i_1,\dots,i_k)$
is a permutation  of $(1,\dots,k)$ and $\d_i$ is a vanishing loop of $\De$ at
the point $e_i,\,\,i=1,\dots,k$. }

\vs

\noin {\it Proof.} Slightly modifying the original homotopy and changing the
numeration of the intersection points $e_1,\,\dots,\,e_n \in {\rm
reg}\,\De$
we may assume that  $e_i \in \gamma_{t_i} \cap \De,\,\,i=1,\dots,k,$
correspond to different
values $0 < t_1 < \dots < t_n < 1$ of the parameter of homotopy
$t \in [0, 1]$.
If $s_i \in [0,\,1],\,\,0 < s_1 < t_1 < \dots < t_n < s_{n+1} < 1,$ and ${\bar
\gamma}_i =
\gamma_{s_i}\,:\,S^1 \to E \se \De , \, i=1,\dots,n+1$, then clearly
${\bar \gamma}_{i+1}^{-1}\cdot {\bar \gamma}_i \approx \d_i$, i.e.
${\bar \gamma}_i \approx
{\bar \gamma}_{i+1}\cdot\d_i$ in $E \se \De$,
where $\d_i$ is a
vanishing
loop of $\De$ at the point $e_i$, and ${\bar\gamma}_1 \approx
\g_0,\,\,{\bar\gamma}_{n+1} \approx \g_1$. Thus, $\gamma_0  \approx
\g_1
\d_n\cdot \dots \cdot \d_1$ in
$E \se \De$, and the lemma follows.  \qed

\vs

In the proof of Theorem 0.2 below  we use the following proposition.
Actually, it follows from Lemma 1.5(C) in [No]. However, we give a
proof which is different from that in [No].

\vs

\noin {\bf 4.3. Proposition.} {\it Let a morphism $p\,:\,E \to B$ of smooth
quasiprojective varieties be a smooth fibration over $B$ with
a connected generic fibre $F$ of positive dimension.
Let $\De \s E$ be a Zariski closed
hypersurface which contains no entire fibre of $p$, i.e.
$p^{-1}(b) \not\subset \De$ for each $b \in B$. Then we have
the following exact sequence}:
$$\pi_1(F \se \De) \stackrel{i_*}{\rightarrow}   \pi_1(E \se \De)
\stackrel{p_*}{\rightarrow} \pi_1(B) \to {\bf 1}\,.$$

\noin {\it Proof.} By Lemma 4.1,
there exist hypersurfaces  $A \s B$ and  $D = H \cup \De \s E$, where
$H := f^{-1}(A)$,    such that   $p\,|\,(E \se D) \,:\, E \se D \to B \se A$
is a
smooth fibration. In particular, $p\,|\,(E \se D)$ induces  an epimorphism of
the fundamental groups. Since the same is also true for
the embedding $i\,:\,B \setminus A \hookrightarrow B$, and since
$p_* = i_* \circ (p\,|\,E \se D)_*$, the exactness at the third term follows.
It remains to prove that the homomorphism
$$i_*\,:\,\pi_1(F \setminus \De) \to {\rm Ker}\,p_* \s \pi_1(E \se \De)$$
is surjective.

Fix a generic fibre $F \not\s D$ and base points
$e_0 \in F \se D$ and $b_0 = p(e_0) \in B \se A$. Let a
class $[\g_0] \in {\rm Ker}\,p_*$ be
represented by a loop $\g_0\,:\,S^1 \to E \se \De$ with the base point
$e_0$. We will show that $\g_0$ is homotopic in $E \se \De$ to a loop
$\g'_0\,:\,S^1 \to F \se \De$ with the same base point.

The loop ${\bar \g}_0 := p \circ \g_0\,:\,S^1 \to B$ with the base point
$b_0 \in B$ is contractible. Let ${\bar \g} \,:\,S^1 \times [0;\,1] \to B$ be
a contraction to the constant loop ${\bar \g}_1 \equiv b_0$. Since
$p\,:\,E \to B$ is a fibration, there exists a covering homotopy
$\g\,:\,S^1 \times [0;\,1] \to E$. Thus, we have
${\bar \g} = p \circ \g$ and
$\g_1\,:\,S^1 \to F$.

Fix a stratification of $D = \De \cup H$ which satisfies the
Whitney condition A and
contains the regular part ${\rm reg}\,D$ of $D$ as an open stratum.
By Thom's Transversality Theorem, the homotopy $\g$ can be chosen
being transversal to the strata of this stratification, and therefore
such that its image
meets the divisor $D$ only in a finite number of its regular points. Let it
meet $\De$ at the points $e_1,\dots,e_k \in$ reg$\,(\De \se H)$. We may also
assume that the loop $\g_1\,:\,S^1 \to F$ does not meet $D$; in
particular, $[\g_1] \in \pi_1 (F \se \De;\,e_0)$. By Lemma 4.2,  $\g_0$ is
homotopic in
$E \se \De$ to the product $\g_1\d_{i_1}\dots \d_{i_k}$, where
$\d_{i}$ is a vanishing loop of $\De$ at the point $e_i,\,\, i=1,\dots,k$.

Note that all the transversal discs to $\De$ in $E$ centered at $e_i$ are
homotopic (via the family of such discs). Hence, all the simple
positive local vanishing loops of $\De$ at $e_i$ are freely homotopic in
$E \se \De$. Therefore, performing further deformation of the vanishing
loops $\d_{i},\,\, i=1,\dots,k$, and taking into account our assumptions
that dim$\,F > 0$ and $\De$ does not contain entirely a fibre of $p$,
we may suppose that

\ss

\noin (i) for each $ \,\, i=1,\dots,k$ the fibre of $p$ through the point
$e_i$ is transversal to $\De$;

\ss

\noin (ii) the loops $\d_{i},\,\, i=1,\dots,k$, do not meet $H$, and
the corresponding local loops $\beta_i,\,\,i=1, \dots, k$, are contained
in the fibres of $p$.

\ss

\noin Since $p^{-1}(A) \s D$, we have that for each
$i = 1,\dots,k$ the projection ${\bar \d}_i := p \circ \d_i$ of the
loop $\d_i$ does not
meet the hypersurface $A \s B$. By the construction, the loops ${\bar
\d}_i,\,\,i = 1,\dots,k,$ are
contractible in $B \se A$. Applying the covering homotopy theorem
to the  smooth fibration $p\,:\,E \se D \to B \se A$ we may conclude that for
each $i = 1,\dots,k,$ the loop $\d_i$ is homotopic in
$E \se D \s E \se \De$ to a loop $\d_i'\,:\,S^1 \to F \se D \s F \se \De$.
Hence, $\g_0$ is homotopic in $E \se \De$ to the product $\g'_0 :=
\g_1\d'_{i_1}\dots \d'_{i_k}\,,\,\,\g'_0\,:\,S^1 \to F \se \De$, and we are
done. \qed

\vs

\noin {\bf 4.4. Duality, discriminants and the Zariski embedding}

\ms

The following construction was used, for instance, in [Zar, pp.307, 326]
and in [DoLib, sect.1, 3]. Let $M$ be an irreducible smooth projective
variety, and let $L \s H^0(M,\,{\cal L})$ be a linear system of effective
divisors on $M$, where $\cal L$ is a linear bundle on $M$. It
defines a rational mapping $\Phi_L\,:\,M \to \pr (L^*)$.
If $K \s L$ is a linear subsystem, then the mapping
$\Phi_K\,:\,M \to \pr (K^*)$ is composed of the mapping $\Phi_L$ followed by
the linear projection $\pi_{L,\,K}\,:\,\pr (L^*) \to  \pr (K^*)$ which is dual
to the tautological embedding $\rho_{K,\,L}\,:\,\pr (K) \hookrightarrow \pr
(L)$.

Set $C_L := \Phi_L (M) \s \pr (L^*)$ and $C_K := \Phi_K (M)
\s \pr (K^*)$, so that $C_K = \pi_{L,\,K} (C_L)$.  The dual variety
$\De_L \s \pr(L)$
of $C_L \s \pr (L^*)$ is usually a hypersurface, which is called {\it the
discriminant hypersurface of the linear system $L$}. The  embedding
$\rho_{K,\,L}$ yields the embedding of the discriminants
$\De_K = \pr (K) \cap \De_L \hookrightarrow \De_L$.

\ss

In particular, starting with a degree $d$ irreducible plane curve $C \s \pr^2$
with a normalization $M \to C$, denote by $K = g^2_d$ the linear system on $M$
of line cuts of $C$ and by $L = |g^2_d|$ the corresponding complete linear
system. Since $g^2_d$ and therefore, also $L$ are base point free, they define
morphisms $\Phi_K\,:\,M \to C \s \pr^2 = \pr (K^*)$ resp.
$\Phi_L\,:\,M \to C_L := \Phi_L(M) \hookrightarrow \pr (L^*)$, and
$C \s \pr^2$ is a projection of the curve $C_L \s \pr (L^*)$. Set $\pr^2_C :=
\pr(K) \hookrightarrow \pr(L)$. The discriminant
$\De_L = C_L^*$ is, indeed, a projective hypersurface, and the dual curve
$C^* \s \pr^2_C$ is an irreducible component of the plane cut $\De_K = \pr(K)
\cap \De_L$. The other irreducible components of $\De_K$ are special tangent
lines of $C^*$ dual to the cusps of $C$ (by {\it a cusp} we mean here a
singular point of a local irreducible analytic branch of $C$). We call these
tangent lines {\it artifacts} [DeZa1]. Thus, the plane cut $\pr(K) \cap \De_L$
of the discriminant hypersurface $\De_L$ is irreducible iff $C$ is an
immersed curve.

The embedding $\pr^{2*} \cong \pr^2_C \hookrightarrow \pr(L)$ which represents
$C^*$ as a plane cut of the discriminant hypersurface $\De_L$ is called
{\it the Zariski embedding} (see [Zar, pp.307, 326; DeZa1]).

By definition, the dual variety $C_L^* = \De_L$ consists of the points
$x \in \pr(L)$ such that the dual hyperplane $x^* \s \pr_L$ cuts out of $C_L$
a non-reduced divisor on the normalization $M$ of $C_L$. If $x \in C_L$ is
a cusp, then, clearly, the dual hyperplane $x^*$ is an irreducible component
of $\De_L$. Thus, the discriminant $\De_L$ is irreducible iff $C_L$ was an
immersed curve. In particular, this is the case if $C =
\pi_{L,\,K} (C_L)$ is an immersed curve. Vice versa, if $C_L$ is an immersed
curve, then the same is true for its generic projection onto the plane. Or,
what is the same, if the discriminant $\De_L$ is irreducible, then its generic
plane section is irreducible, too.

The projectivization $\pr(L)$ of the complete linear system $L$ of degree $d$
divisors on $M$ coincides with a fibre of the Abel--Jacobi map
$\phi_d\,:\,S^dM \to J_d(M)$ (see (3.2)), so that $\pr^2_C$ is a plane in this
fibre. We still call the morphism $\pr^2_C \hookrightarrow S^dM$
{\it the Zariski embedding}. The hypersurface $\De_d \s S^dM$ which consists of
the non-reduced degree $d$ effective divisors on $M$ is also called
{\it the discriminant hypersurface}. It is the image of the diagonal
hypersurfaces of the direct product $M^d$ via the Vieta map $M^d \to S^dM$.
Thus, $\De_L = \pr(L) \cap \De_d$, where $\pr(L)$ has been identified with a
fibre $F_j := \phi_d^{-1}(j),\,j \in J_d(M)$, of $\phi_d$.

\vs

\noin {\bf 4.5.} {\it Proof of Theorem 0.2.} By Corollary 3.4, we may
suppose that $C$ is a generic nodal Pl\"ucker curve of degree $d$ and
geometric genus $g$, where $d \ge 2g-1$. By Mattuck's Theorem  (see (3.2)),
the $d$-th Picard bundle
$\phi_d\,:\,S^dM \to J_d(M)$, where $M$ is a normalization of $C$, is a
projective bundle with a generic fibre
$F \cong \pr^{d - g}$. By (4.4), the dual curve $C^*$ can be identified
with the plane cut of the discriminant hypersurface $\De_d \s S^dM$
by the plane $\pr^2_C$ via its Zariski embedding $\pr^2_C \hookrightarrow F_0
:= \pr(L) \s S^dM$, where $L = |g^2_d|$ and  $g^2_d$ is the linear system
on $M$ of line cuts of $C$.

If the group $\pi_1 (\pr^2_C \se \Delta_d)$ is big for a generic plane
$\pr^2_C \s F_0 \cong \pr^{d - g}$, then by Zariski's Lefschetz type Theorem
[Zar, p.279; Di, 4.1.17], it is big for any such plane, so that
$\pr^2_C \s F_0$ might be assumed being generic. Indeed, a section $S$ of the
discriminant hypersurface $\De_L = F_0 \cap \De_d$ by a generic plane $\pr \s
F_0$ is an irreducible curve with the same normalization $M$ and with the dual
$S^* \s \pr^2$ an immersed curve of degree $d$ and genus $g$. Thus, we may
start with $C = S^*$ and obtain $C^* = S = \pr^2_C \cap  \De_d$. Note that
such a generic linear system $K = g^2_d \s L$, where $\pr = \pr(K)$, defines
a morphism $M \to \pr^2$ such that its image coincides with $C = S^*$.
Since by Theorem 2.1(c), $PlNod_{d,\,g}$ is a Zariski open subset of
$Imm_{d,\,g}$, the curve $C=S^*$ obtained in this way is a nodal Pl\"ucker one.

Applying the Zariski Lefschetz type Theorem we get an isomorphism
$$\pi_1 (F_0 \se \Delta_d) \cong \pi_1 (\pr^2_C \se \Delta_d) \cong
\pi_1 (\pr^2 \se C^*)\,.$$
By Proposition 4.3, we have
the exact sequence $$\pi_1 (F_0 \se \Delta_d) \to \pi_1 (S^d M \se  \Delta_d)
\to \pi_1 (J_d(M)) \cong \gz^{2g} \to {\bf 1}\,.$$
It follows that $ \pi_1 (\pr^2 \se C^*)$ is a big group if $\pi_1 (S^d M \se
\Delta_d)$ is big (cf. (1.1)). But $\pi_1 (S^d M \se  \Delta_d)$
is the braid group $B_{d,\,g}$ of $M$ with $d$ strings which is big (see
Lemma 1.2($b$)). This completes the proof. \qed

\vs

\noin {\it Remark.} A presentation of the group $\pi_1(\pr^2 \se C)$
for a generic maximal cuspidal curve $C \s \pr^2$ of genus  0 or 1
was found by Zariski [Zar, p. 307]; see also [Ka] for $g \le {d-1 \over 2}$,
where $d = {\rm deg}\,C^*$. The result of [Ka] is based on the statement in
[DoLib] that for $d \ge 2g-1$ the $d$--th
Abel--Jacobi mapping $\phi_d\,:\,S^dM \to J(M)$ restricted to the
complement of the discriminant hypersurface $\De_d \s S^dM$ is a Serre
fibration, so that  the long exact homotopy sequence is available. But the
indication given in [DoLib] does not seem to be sufficient
for the proof. Another proof of the exactness of the above sequence of
fundamental groups extended to the left by the term $\bf 1$ has been recently
obtained in [KuShi].  Once again, this leads to a presentation of
the group $\pi_1(\pr^2 \se C)$.

\vs

\begin{center} {\LARGE References} \end{center}

{\footnotesize

\noin [ACGH] E. Arbarello, M. Cornalba, P.A. Griffiths, J. Harris, {\sl
Geometry of algebraic curves. I}, N.Y. e.a.: Springer, 1985

\noin [AC] E. Arbarello, M. Cornalba. {\sl Su una conjettura di Petri},
Comment. Math. Helv. 56 (1981), 1-38

\noin [Au] A. B. Aure. {\sl Pl\"ucker conditions on plane rational curves},
Math. Scand. 55 (1984), 47--58, with {\sl Appendix} by S. A. Str\"omme,
ibid. 59--61

\noin [BPVV] W. Barth, C. Peters, A. Van de Ven. {\sl Compact complex
surfaces},
N.Y. e.a.: Springer, 1984

\noin [Be] D. Bennequin. {\sl Entrelacements et \'equations de Pfaff},
Ast\'erisque, 107--108 (1983), 87--161

\noin [BinFl] J. Bingener, H. Flenner. {\sl On the fibres of analytic
mappings}, in: Complex Analysis and Geometry, V. Ancona and A. Silva, eds.,
N.Y.: Plenum Press, 1993, 45--101

\noin [Bi] J.S. Birman, {\sl Braids, links, and mapping class groups},
Princeton Univ. Press, Princeton, NJ, 1974

\noin [Bo] A. Borel. {\sl K\"ahlerian coset spaces of semisimple Lie groups},
Proc. Nat. Acad. Sci., USA, 40, No. 12 (1954), 1147--1154

\noin [BoHC] A. Borel, Harish--Chandra. {\sl Arithmetic subgroups of
algebraic groups}, Ann. of Math. 75 (1962), 485--535

\noin [CoZi] D. J. Collins, H. Zieschang. {\sl Combinatorial group theory and
fundamental groups}, In: Algebra VII, Encyclopaedia of Math. Sci. Vol. 58,
Berlin e.a.: Springer, 1993, 3--166

\noin [DeZa1] G. Dethloff, M. Zaidenberg. {\sl Plane curves
with hyperbolic and C--hyperbolic complements},
Ann. Sci. Ecole Norm. Super. Pisa (to appear); Pr\'epublication de l'Institut
Fourier de Math\'ematiques, 299, Grenoble 1995, 44p.

\noin [DeZa2] G. Dethloff, M. Zaidenberg. {\sl  Examples of plane curves of low
degrees with hyperbolic and C--hyperbolic complements.}  In:
Geometric Complex Analysis, J. Noguchi e.a. eds. World Scientific Publ. Co.,
Singapore 1996., 176--193

\noin [DeOrZa] G. Dethloff, S. Orevkov, M. Zaidenberg.
{\sl Plane curves with a
big fundamental group of the  complement.} Pr\'epublication de l'Institut
Fourier de Math\'ematiques, 354, Grenoble 1996, 26p.; E-print alg-geom/9607006

\noin [Di] A. Dimca, {\sl Singularities and topology of hypersurfaces},
Berlin e.a.: Springer, 1992

\noin [DoLib] I. Dolgachev, A. Libgober. {\sl On the fundamental group of the
complement to a discriminant variety}, In: Algebraic Geometry,
Lecture Notes in Math. 862, 1--25, N.Y. e.a.: Springer, 1981

\noin [GoShVi] V. V. Gorbatsevich, O. V. Shvartsman, E. B. Vinberg.
{\sl Discrete subgroups of Lie groups.} In: Lie Groups and Lie Algebras II,
Encyclopaedia of Math. Sci. Vol. 21, Springer : NY e.a., 1995

\noin [Ha] J. Harris. {\sl On the Severi problem}, Invent. Math. 84 (1986),
445--461

\noin [Hart] R. Hartshorn. {\sl Algebraic geometry}, NY e.a.: Springer, 1977

\noin [He] S. Helgason. {\sl Differential geometry, Lie groups and
symmetric spaces}, N.Y. e.a.: Academic Press, 1978

\noin [Iv] N. V. Ivanov. {\sl Algebraic properties of the Teichm\"uller
modular group}, Soviet Math. Dokl. 29 (1984), No.2, 288--291

\noin [Ka] J. Kaneko. {\sl On the fundamental group of the complement to
a maximal cuspidal plane curve}, Mem. Fac. Sci. Kyushu Univ. Ser. A. 39
(1985), 133-146

\noin [Ko] K. Kodaira. {\sl A theorem of completeness of characteristic
systems for analytic families of compact submanifolds of complex manifolds},
Ann. Math. 75 (1962), 146--162

\noin [Kos] J.-L. Koszul. {\sl Sur la forme hermitienne canonique des espaces
homog\`enes complexes}, Can. J. Math. 7 (1955), 562--576

\noin [KuShi] Vik. S. Kulikov, I. Shimada. {\sl On the fundamental group of
complements to the dual hypersurfaces of projective curves.} Preprint
Max-Planck-Institute f\"ur Mathematik, MPI 96-32, Bonn, 1996,  1--15

\noin [Lib] A. Libgober. {\sl Fundamental groups of the complements to plane
singular curves}, Proc. Sympos. in Pure Mathem. 46 (1987), 29--45

\noin [Lin] V. Ja. Lin. {\sl Liouville coverings of complex spaces, and
amenable groups}, Math. USSR Sbornik, 60 (1988), 197--216

\noin [LinZa] V. Ja. Lin,  M. Zaidenberg. {\sl Liouville and Carath\'eodory
coverings in Riemannian and complex geometry.} Preprint Max-Planck-Institute
f\"ur Mathematik, MPI 96-110, Bonn, 1996, 1--20 (see in this volume)

\noin [LySu] T. Lyons, D. Sullivan. {\sl Function theory, random paths,
and covering spaces}, J. Diff. Geom. 19 (1984), 299--323

\noin [Ma] A. Mattuck. {\sl Picard bundles}, Illinois J. Math. 5 (1961),
550--564

\noin [MC] J. McCarthy. {\sl A "Tits-alternative" for subgroups of surface
mapping class groups}, Trans. Amer. Math. Soc. 291 (1985), 583--612

\noin [Mi] J. Milnor. {\sl Singular points of complex hypersurfaces},
Princeton, New Jersey : Princeton University Press, 1968

\noin [MoTe] B. Moishezon, M. Teicher. {\sl Fundamental groups of complements
of branch curves as solvable groups}, Duke E-print alg--geom/9502015, 1995,
17p.

\noin [Na] M. Namba. {\sl Geometry of projective algebraic curves},
N.Y. a.e.: Marcel Dekker, 1984

\noin [No] M.V. Nori. {\sl Zariski's conjecture and related problems},
Ann. scient. Ec. Norm. Sup. 16 (1983), 305--344

\noin [O] M. Oka. {\sl Symmetric plane curves with nodes and cusps}, J.
Math. Soc. Japan, 44, No. 3 (1992), 375--414

\noin [OSh] A. Yu. Ol'shanskij, A. L. Shmel'kin. {\sl Infinite groups},
In: Algebra IV, Enciclopaedia of Math. Sci. 37, Berlin e.a.: Springer, 1993,
3--95

\noin [Ra] M. S. Raghunathan. {\sl Discrete subgroups of Lie groups},
Berlin e.a.: Springer, 1972

\noin [Se] F. Severi. {\sl Vorlesungen \"{u}ber algebraische Geometrie},
Leipzig: Teubner, 1921

\noin [Sh1] G. B. Shabat. {\sl The complex structure of domains covering
algebraic surfaces}, Functional Analysis Appl. 11 (1977), 135--142

\noin [Sh2] G. B. Shabat. {\sl On families of curves covered by bounded
symmetric domains}, Serdica Bulg. Math. Public. 11 (1985), 185--188
(in Russian)

\noin [Ti] J. Tits. {\sl Free subgroups in linear groups}, J. Algebra, 20
(1972), 250--270

\noin [Va] A.N. Varchenko. {\sl Theorems on the topological equisingularity
of families of algebraic varieties and families of polynomial mappings},
Math. USSR Izvestija, Vol. 6 (1972), No. 5, 957--1019

\noin [Zar] O. Zariski. {\sl Collected Papers}. Vol III : {\sl Topology of
curves and surfaces, and special topics in the theory of algebraic
varieties},  Cambridge, Massachusets e. a.: The MIT Press, 1978

\vs

\noindent Gerd Dethloff,
Mathematisches Institut der Universit\"at G\"ot\-tin\-gen,
Bunsenstrasse 3-5,
37073 G\"ot\-tin\-gen,
Germany.
e-mail: DETHLOFF@CFGAUSS.UNI-MATH.GWDG.DE

\vs

\noin Stepan Orevkov,
System Research Institute RAN,
Moscow, Avtozavodskaja 23, Russia.
e-mail: OREVKOV@MAIN.DOMINO.MSK.RU

\vs

\noin  Mikhail Zaidenberg,
Universit\'{e} Grenoble I,
Institut Fourier et Laboratoire de Math\'ematiques
associ\'e au CNRS,
BP 74,
38402 St. Martin d'H\`{e}res--c\'edex,
France.
e-mail: ZAIDENBE@PUCCINI.UJF--GRENOBLE.FR}

\end{document}